\pdfoutput=1
\documentclass{iopart}
\usepackage{amsthm, amssymb}
\usepackage{graphicx}

\usepackage{mathptmx}

\usepackage{subfigure}

\renewcommand{\phi}{\varphi}
\renewcommand{\leq}{\leqslant}
\renewcommand{\geq}{\geqslant}
\newcommand{\defeq}{:=}

\newcommand{\norm}[1]{\left \| #1 \right \|}
\newcommand{\nn}{\mathbf{n}}
\newcommand{\qq}{\mathbf{q}}
\newcommand{\vv}{\mathbf{v}}
\newcommand{\ttt}{\mathbf{t}}
\newcommand{\xx}{\mathbf{x}}

\newcommand{\figref}[1]{figure~\ref{#1}}
\newcommand{\bfigref}[1]{Figure~\ref{#1}}
\newcommand{\eqref}[1]{(\ref{#1})}

\usepackage{url}
\begin{document}

\vspace*{-50pt}
\title[Bifurcations of periodic and chaotic attractors in pinball billiards]{Bifurcations of periodic and chaotic attractors in pinball billiards with focusing boundaries}

\author{Aubin Arroyo$^1$, Roberto Markarian$^2$ and David P.~Sanders$^3$}
\address{$^1$  Instituto de Matem\'aticas, Unidad Cuernavaca, Universidad Nacional Aut\'onoma de M\'exico,
Apartado postal 273-3, Admon.\ 3, Cuernavaca, 62251 Morelos, Mexico}
\address{$^2$ Instituto de Matem\'atica y Estad\'istica (IMERL), Facultad de Ingenier\'ia, Universidad de la Rep\'ublica,
Montevideo 11300, Uruguay}
\address{$^3$ Departamento de F\'isica, Facultad de Ciencias, Universidad Nacional Aut\'onoma de M\'exico, Ciudad Universitaria, 04510 M\'exico D.F., Mexico}
\ead{aubin@matcuer.unam.mx\textrm{,} roma@fing.edu.uy \textrm{and}
dps@fciencias.unam.mx}

\begin{abstract}
We study the dynamics of billiard models with a modified collision rule: the outgoing angle from a collision is a uniform contraction, by a factor $\lambda$, of the incident angle.
These \emph{pinball billiards} interpolate between a one-dimensional map when $\lambda=0$ and the classical Hamiltonian case of elastic collisions when $\lambda=1$. For all $\lambda<1$, the dynamics is dissipative, and thus gives rise to attractors, which may be periodic or chaotic. 
Motivated by recent rigorous results of Markarian, Pujals and Sambarino \cite{MarkarianPujalsSambarinoPinballBilliards},
we numerically investigate and characterise the bifurcations
of the resulting attractors as the contraction parameter is varied. Some billiards exhibit only periodic attractors, some only chaotic attractors, and others have coexistence of the two types.
\end{abstract}



\section{Introduction}

Billiard models have played an important role in mathematical
physics for forty years, since their introduction by Sinai in
\cite{SinaiElasticReflectionsRussMathSurv1970}, motivated by the Boltzmann
ergodic hypothesis on fluids of hard spheres.
They are a class of models which are accessible and of interest both to mathematicians \cite{ChernovMarkarianChaoticBilliardsAMS2006} and to physicists \cite{GaspardBook1998}.
In particular, they have attracted much attention as simple models of physical
systems which exhibit strong chaotic properties (hyperbolicity, ergodicity and
mixing)

In classical billiards, particles travel freely between collisions with hard, fixed boundaries. They then collide \emph{elastically}, i.e.\ without friction, so that the kinetic energy is conserved. Such billiards can thus be considered as a type of Hamiltonian system with discontinuities.

Our objective is to investigate, mainly numerically, the different types of dynamics
which occur in billiard tables with a \emph{modified} collision
rule: we consider \emph{pinball billiards}, where the reflection
rule at collisions is no longer elastic. Rather, the angle (with respect
to the normal) of reflection is smaller than the angle of
incidence, by a constant factor. The dynamics is then no longer conservative; rather,
 there is contraction in phase space, which gives rise to attractors.
Two types of attractors are found to occur: they are either (i) periodic points
whose stability matrix has eigenvalues with modulus not greater
than one, or (ii) non-trivial invariant sets that exhibit 
\emph{dominated splitting} (defined below). Both types of attractors can coexist,
and we study examples where all possible cases occur.

Rather than attempt to model a particular physical scenario, we consider a model that interpolates between a one-dimensional dynamics, where reflection is always normal to the boundary, and a Hamiltonian one (elastic reflection).
A related modified billiard collision rule, which was motivated by physical considerations, was studied recently in 
\cite{AltmannDelMagnoWaveInspiredCorrectionsBiliardsEPL2008}.
Other Hamiltonian systems with dissipation have been studied from a similar point of view, for example the kicked rotor map \cite{FeudelGrebogiMap100CoexistingPeriodicAttractorsPRE1996, GallasMultistabilityKickedRotorIJBC2008}.

We are interested in the dynamical properties of orbits in the forward limit set of the system, that is, the union of all $\omega$-limit sets of points in the space, which describes the dynamics of the system as the time goes to $+\infty$. Markarian, Pujals and Sambarino \cite{MarkarianPujalsSambarinoPinballBilliards} proved that, under certain circumstances, case (ii), that is the ocurrence of non-trivial invariant sets with dominated splitting, appears in many common billiard tables. In our context, an invariant subset $\Lambda$ of the limit set has \emph{dominated splitting} if the tangent bundle continuously splits into two non-trivial invariant sub-bundles under the derivative, say $T_{\Lambda} = E \oplus F$, and there are $C>0$ and $\nu \in (0,1)$ such that
\begin{equation}
\norm{Df^{n} |_{E(x)} } \cdot \norm{Df^{-n} |_{F(f^n(x))}} \leq C \nu^n
\end{equation}
for any $x \in \Lambda$ and $n\geq 0$. In these circumstances one can expect \cite{MarkarianPujalsSambarinoPinballBilliards} that the limit set can be decomposed into the union of three subsets: (1) a finite union of disjoint invariant curves formed by periodic points with finite period; (2) a finite union of simple closed curves where the dynamic is conjugate to an irrational rotation; and (3) a disjoint union of finitely many compact invariant hyperbolic sets.

In this paper we restrict attention to tables with \emph{focusing}
boundaries. A focusing boundary (which has negative curvature under our convention) has the property that a beam of parallel trajectories which hits it focuses after the collision.
Bunimovich showed \cite{BunimovichErgPropsNowhereDispersingBilliardsCMP1979} that
if there is ``enough space'' in the billiard table, then after it has focused, the beam can then defocus, before it collides again. This \emph{defocusing mechanism} 
can give rise to strong chaotic properties, in particular hyperbolicity \cite{WojtkowskiPrinciplesDesignBilliardsLyapExpsCMP1986, MarkarianBilliardsPesinRegionMeasureOneCMP1988} and to systems which are ergodic and mixing, see for instance \cite{SzaszKPropertyPlanarHypBilliardsCMP1992} and \cite{MarkarianNewErgodicBilliardsCardioidNonlin1993}.
 
 One of the motivations of this work is to study pinball billiard maps defined on tables which in some sense are close to the limit of validity of the hypotheses assumed in \cite{MarkarianPujalsSambarinoPinballBilliards} for which rigorous results could be obtained. We consider cases in which either the attractor is not contained in a compact region of the phase space, or certain geometrical properties used in that reference do not hold.  In this way,  we shed light on the interesting types of behaviour which pinball billiards can exhibit.
 
We begin by defining formally the dynamical system, and then proceed to study numerically the dynamics in several different types of billiard table with
focusing (and in one case neutral) boundaries.
The aim of the paper is to characterise the types of attractor that can be expected to occur in pinball billiards,
and relate them to the geometric features of the billiard table. In particular, there is a competition between the defocusing mechanism and the angle contraction
in the pinball dynamics.

  The tables we consider are the family of ellipses, the cardioid, a modified, convex version of the cardioid with a discontinuity in the curvature, and a smooth family of perturbations of the circle.
In the ellipse, it is found that the limit set is always periodic and it can either have dominated splitting or not, depending on whether certain eigenvalues are real or complex.
   In the cardioid, the limit set always consists of a chaotic attractor which accumulates on the vertex.
   We then introduce a modified, convex version of the cardioid, in which the vertex is removed. There we find that the limit set consists of periodic or chaotic attractors, and for some values of the contraction parameter these two phenomena coexist.   Here there are two points where the curvature is discontinuous, and the chaotic attractor contains trajectories which include these points. 
   Finally, in a smooth family of perturbations of the circle we observe coexistence of periodic and chaotic attractors, which undergo bifurcations when the contraction parameter is varied.

\section{Pinball billiards}
\label{sec:definition}

Let $Q$ be an open, bounded, simply-connected subset of the plane,
whose boundary $\Gamma$ consists of a finite number of compact $C^3$ curves $\Gamma_i$. 
\emph{Billiard} dynamics in $Q$ is the dynamical system describing the free motion of a point mass particle inside $Q$ with reflections at the boundary $\Gamma \defeq \cup_i \Gamma_i$.
Each $\Gamma_i$ is called a \emph{(smooth) component} of $\Gamma$. Let $n(q)$ be the unit normal vector at the point $q\in \Gamma_i$ which points towards the
interior of $Q$.  The \emph{configuration space} of this dynamical system is $Q$, and its  \emph{phase space} $M$ is given by
\begin{equation}
 M \defeq \{(q, v): q \in \Gamma, |v|=1, \langle v, n(q) \rangle \geq 0 \},
\end{equation}
where $\langle \cdot \rangle$ denotes the Euclidean inner product.
 
We introduce coordinates $(s, \phi)$ on the phase space $M$, where $s$
is the arc length parameter along $\Gamma$ and $\phi$ is the angle
between the reflected vector $v$ and the normal $n(q)$, satisfying
$\langle v, n(q) \rangle = \cos\phi$. In classical billiards the
angle of incidence is equal to the angle of reflection (elastic
reflection) \cite{ChernovMarkarianChaoticBilliardsAMS2006}; in this case, the
billiard map is area-preserving in terms of the coordinates $(s, \sin(\phi))$,
known as Birkhoff coordinates.

Recently, Markarian, Pujals and Sambarino \cite{MarkarianPujalsSambarinoPinballBilliards}
introduced \emph{pinball billiards}, which are 
billiard maps with a modified reflection rule. Let $\eta \in
[-\frac{\pi}{2}, \frac{\pi}{2}]$ be the angle from the incidence
vector to the outward normal $-n(q)$ at the point $q$ where the
ball hits the boundary. In pinball billiards we no longer restrict to the elastic case; rather,
we allow the outgoing angle $\phi$ to depend on the incidence angle $\eta$. In principle,
a dependence could be allowed on the position of the collision, as in \cite{AltmannDelMagnoWaveInspiredCorrectionsBiliardsEPL2008}, but 
we suppose that the
exit angle $\phi$ depends only on the incidence angle $\eta$. Then $\phi =
-\eta + f(\eta)$, where $f: [-\pi/2, \pi/2] \to \mathbb{R}$ is a
$C^2$ function.

In \cite{MarkarianPujalsSambarinoPinballBilliards}, section~3.2 of which gives more details on
these topics, two perturbations on the angle of reflection were
studied. In this paper, we restrict attention to one of them: we assume
that $f(\eta) = (1 - \lambda)\eta$, with constant $\lambda <1$; thus
$\phi = -\lambda \eta$, 
that is, the outgoing angle is given by a \emph{uniform} contraction of 
the angle of incidence.  For simplicity we rescale the speed of the particle after each collision so that it is always $1$.
Physically this corresponds to a sticky boundary which gives exactly the right amount of ``kick'' to the particle; this can be thought of as a caricature
of what happens in a real pinball table, whence the term ``pinball billiard''.

The \emph{pinball billiard map} $T_{\lambda}$ with a given contraction parameter $\lambda \in [0,1]$ is then defined as follows: we leave the
point $q_0 \in \Gamma$ with direction $v_0$ until we intersect the
boundary again at the point $q_1 \in \Gamma$. If the unit normal $n(q_1)$ 
at the new position is well defined, then the new velocity 
vector $v_1$ is given by the above collision rule. The
pinbilliard map $T_\lambda$ is then defined by $T_\lambda(q_0, v_0) \defeq (q_1, v_1)$. If
$n(q_1)$ is not uniquely defined, for example if the trajectory lands at 
a sharp corner, then the map $T_\lambda$ is not defined at $(q_0, v_0)$.
$T_{\lambda}$ is also undefined at points $(q_{0},v_{0})$ where the trajectory
hits tangentially another component of the boundary $\Gamma$. 
Removing such points from the domain of definition of the map forces us to work in a non-compact domain. This is an obstacle  for the direct application of results from \cite{MarkarianPujalsSambarinoPinballBilliards}.

Let $t$ be the Euclidean distance between $q_0$ and $q_1$, and $K_i$ be
the curvature of $\Gamma$ at $q_i$ taken with the following convention:
$K_i$ is negative if the component of the boundary is focusing,
zero if it is flat, and positive if it is dispersing.  Zero
curvature is not allowed on focusing and dispersing components of
the boundary. 
In this paper we will only study billiards with
focusing and flat components, and only for $\lambda \leq 1$. Although it is possible to consider $\lambda > 1$, in this case the trajectories tend to converge to ``whispering gallery'' modes tangent to the billiard boundary.

\subsection{Derivative matrix}

The total derivative of the pinball billiard map $T_\lambda$ at a point $x=(s, \phi_0) \in M$ is, as shown in \cite[subsection 3.3]{MarkarianPujalsSambarinoPinballBilliards}, given by:
\begin{equation} 
D_{x}T_{\lambda}(s, \phi_0) = - \left(
\begin{array} {cc} A & B
\\ \lambda \left(K_1 A + K_0 \right)   &
\lambda \left ( K_1 B +1 \right)  
\end{array} \right),
\label{eq:derivative_of_map}
\end{equation}
where
 \begin{equation}
  A \defeq \frac{ t K_0+ \cos \phi_0}{ \cos\eta_1}; \qquad B \defeq \frac{t}{\cos \eta_1}.
 \end{equation}
Again, subscript $0$ refers to quantities evaluated at the starting collision point, and subscript $1$ to quantities at the next collision.
 Thus $\phi_0$ is the exit angle of the trajectory at the boundary point from which the trajectory leaves, and $\eta_1$ is the incidence angle at the \emph{next} collision, calculated by the rule described above.  The derivative thus involves information from two consecutive collisions.
 
 The stability matrix of an orbit is the product of the stability matrices of each trajectory segment. 
In the case of a periodic orbit, it is sufficient to take the product over a single period.
 
\subsection{Slap billiard map: $\lambda=0$}

The pinball billiard map depends on the contraction parameter 
$\lambda$, and one of our goals is to study the global modifications (bifurcations)
of the dynamics of $T_\lambda$ as $\lambda$ varies in the interval $[0,1]$.
When $\lambda=1$ the outgoing angle is equal to the incidence angle, so that $T_1$ is the classical elastic billiard map, and the measure $\rmd \mu = \cos \phi \, \rmd s \, \rmd \phi$ is preserved.

When $\lambda = 0$, trajectories leave exactly along the normal vector at each collision, so that $T_0$
reduces to a \emph{one-dimensional} map, 
which was named the \emph{slap billiard map} in
\cite{MarkarianPujalsSambarinoPinballBilliards}; the derivative of $T_0$ at a
point $q$ along the boundary
is given by $T_0'(q) = -(t K_0 +1)/ \cos \eta_1$. 
Note that the slap billiard map is in fact a map of the circle, since the arc length is a periodic variable.
We emphasise that in the cases considered in this paper the curvatures $K_i$ are always non-positive. 

We remark that the slap billiard map is  related to the geometrical concept of \emph{evolute} of a curve, which is the envelope of the set of all normals
to the curve. Since the trajectories in the slap billiard always follow the normal vectors, the evolute is generated in configuration space.
Nonetheless, this seems not to be particularly relevant for the dynamical properties of the map.

\section{Circular and elliptical tables}

The simplest smooth, convex tables with focusing boundaries 
are elliptical tables,  which are classical examples in the case of elastic billiards \cite[Chapter~1]{ChernovMarkarianChaoticBilliardsAMS2006}. 
In this section we study pinball billiards in elliptical tables, beginning with the special case of circular tables.
The tables have elliptical boundaries, parametrised as $\Gamma_a(t) =(a \cos t,  \sin t)$, where $t \in [0,2\pi)$ and the 
shape parameter $a$ satisfies $a \geq 1$ without loss of generality. 

\subsection{Circular table}

The simplest case of pinball
billiard dynamics occurs in a circular table, for $a=1$. When
$\lambda = 1$, that is, for elastic collisions, the incident and
outgoing angles are constant along each trajectory, due to conservation of
angular momentum, and the phase
space is foliated by parallel lines corresponding to constant
angles. This system is thus completely regular (integrable).

\begin{figure}
\subfigure[]{ 
\includegraphics*[scale=1]{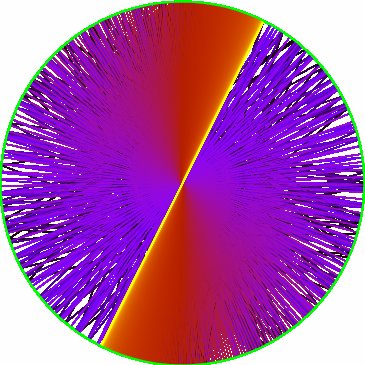}
}
\qquad
\subfigure[]{ 
\includegraphics*[scale=1]{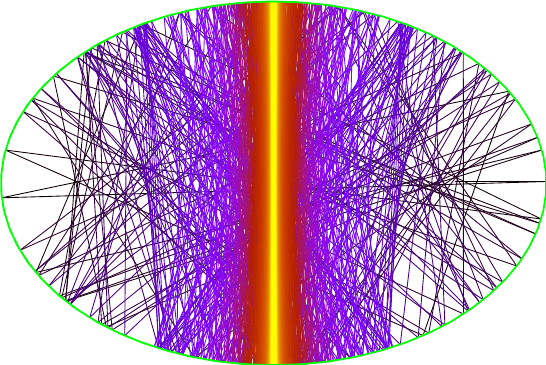}
}
\caption{Single trajectories in configuration space for pinball billiard dynamics with contraction parameter $\lambda=0.99$ in (a) circular table, and (b) elliptical table
with $a=1.5$.
Colours indicate the number of bounces, with
lighter colours corresponding to later times, exhibiting the asymptotic convergence to period-$2$ orbits. 
The initial condition in (a) is a random one; that in (b) was taken close to the unstable period-$2$ orbit along the major axis, from which it rapidly diverges.
Remnants of the structure of caustics which is found for the integrable case ($\lambda=1$) are visible in (b).
}
\label{fig:circle-ellipse-config}
\end{figure}

For $\lambda < 1$, however, at each collision the exit angle
$\phi' = \lambda \phi$ decreases, so that the outgoing
trajectory is closer to the diameter which corresponds to the
normal direction at the point of collision. This continues at each collision; 
since there is no counter-effect, the trajectory from any
initial condition  converges to some period-$2$
diametrical orbit, as shown in \figref{fig:circle-ellipse-config}(a).
Any diameter is a limit set, and the exact limiting behaviour depends on the initial conditions.

In \cite{MarkarianPujalsSambarinoPinballBilliards} it 	is proved that for $0 < \lambda < 1$ the limit set
of the pinball billiard in a circular table has
dominated splitting. This limit set consists of the line
$[\phi=0]$, the set of all period-$2$ orbits.

The slap billiard map of a circular table corresponds to a translation by $\pi$
on the circle: $T_0(s) = s + \pi \mathrm{\ (mod\ } 2 \pi)$, that is, each point $s$ is mapped to 
the diametrically opposed point $s+\pi$ and then back, so that each point is periodic with period $2$.

\subsection{Elliptical table}
New dynamical features appear when we break the circular symmetry and consider 
an ellipse with non-zero eccentricity. An ellipse has exactly two periodic orbits of period
$2$, lying along its major and minor axes. However, these two
periodic orbits are found to have different stability properties under the pinbilliard dynamics: the orbit
lying along the minor axis is \emph{stable}, whereas the orbit
lying along the major axis is \emph{unstable}: an initial
condition pointing along the major axis stays there forever,
whereas an initial condition starting arbitrarily close to it leaves the neighbourhood
of the major axis and eventually converges to the stable orbit, as shown in \figref{fig:circle-ellipse-config}(b). 

The complete description of the dynamics is that there is one stable and one unstable period-2 orbit and they form the limit set of the pinball billiard transformation, for all $\lambda \in (0,1)$. With the exception of the unstable orbit, all other trajectories are attracted by the stable one.

\subsection{Stability of period-$2$ orbits}
To confirm the results on the stability of the period-$2$ orbits, we calculate explicitly their stability matrices.
The curvature of the ellipse boundary at the point $\Gamma_a(t)$ is given by
\begin{equation}
\kappa(t) =\frac{-a}{ \left[ {a}^{2}  \sin^2(t)  +   \cos^2(t)   \right] ^{3/2} }
\end{equation}
Thus the curvature at the points of the ellipse lying on the major
axis is $\kappa(0) = -a$, and on the minor axis is $\kappa(\pi/2) = - 1/a^{2}$.

The stability matrix of the periodic orbit of period $2$ along the minor axis, for
contraction parameter $\lambda$, is the product of the stability matrices for the two flights of the orbit, and is found to be
\begin{equation} \label{monodromyEllipse}
\fl
\left(
 \begin {array}{cc}
  \left[4(1+\lambda)(1-a^2) + a^{4} \right]  / a^{4} \quad &
 2 (1+ \lambda) (a^2-2)  /  a^{2} \\
 \noalign{\medskip}
 -2 \lambda  ( 1+\lambda )  ( a^2 - 1 )   ( a^2 - 2 )/ a^6 \quad &
 \lambda \left[4 (1-a^2)(1+\lambda) + \lambda a^4 \right] / a^{4}
\end{array}
\right).
\end{equation}
The stability of the orbit is determined by the eigenvalues  of \eqref{monodromyEllipse}, which are given by
\begin{equation} \label{eigenvaluesEllipse}
\mu_{\pm}(a,\lambda) = \frac{P_{1}(a,\lambda) \pm \sqrt{P_{2}(a,\lambda)}}{2a^4},
\end{equation}
where
\begin{eqnarray}
P_1  = (a^2-2)^2 \lambda^2 - 8(a^2-1)\lambda + (a^2-2)^2;\nonumber\\
P_2  = [(1 + \lambda) (a^2-2)]^2 [(a^2-2)^2 \lambda^2 - 2(a^4+ 4 a^2 - 4)\lambda + (a^2-2)^2].
\end{eqnarray}
The eigenvalues are complex conjugates satisfying $|\mu_{\pm}(\lambda)| < 1$ for $\lambda$ close to 1 and any $a \geq 1$, so that the minor axis is then a stable focus.
Furthermore, $|\mu_{\pm}(\lambda)| \to 1$ when $\lambda \to 1$, since the orbit
is elliptic in the elastic case $\lambda=1$. 

The eigenvalues change type from complex to real when $P_2(a,\lambda)$ crosses $0$. This occurs when $a = \sqrt{2}$, for any $\lambda \in [0,1]$.
For $a \neq \sqrt{2}$, we have also the following solutions:
\begin{equation}
\lambda_{\pm}(a) =  \frac{a^4 +4a^2-4 \pm 4 a^2 \sqrt{a^2-1}}{(a^2-2)^2}.
\end{equation}
We find that $\lambda_{+}$ is always outside the region of $\lambda$ of interest, but $\lambda_{-}(a) \in [0,1]$ for all $a\geq 1$, with $\lambda_{-}(a) \to 0$ when $a \to \sqrt{2}$ and $\lambda_{-}(a) \to 1$ when $a \to \infty$. For pairs $(a,\lambda)$ below the graph of $\lambda_{-}(a)$, both $\mu_{\pm}(a,\lambda)$ are real and the dynamics is then that of a stable node; above this graph, both eigenvalues have non-trivial imaginary part.  These eigenvalues are shown for the case $a=5$ in \figref{fig:ellipse-evals}.
Note that the absolute value of the eigenvalues must reach $1$ when $\lambda=1$, since the orbit is elliptic in the elastic case.


The stability of the other period-$2$ orbit, along the major axis, can be found by inverting the geometry.
A computation shows that both eigenvalues $\mu_\pm$ are real, with $0<\mu_- < 1 < \mu_+$. Hence the major axis is always an unstable period-$2$ orbit for $\lambda<1$.
 
In our two-dimensional setting, the first requirement to have a set with dominated splitting is the existence of two different non-trivial invariant directions for $DT_{\lambda}$. Thus, if the limit set contains a periodic point with non-real eigenvalues, then this set cannot have dominated splitting. On the other hand, if the limit set consists of only a finite number of hyperbolic periodic orbits with real eigenvalues, then the limit set does have dominated splitting.
Hence, for any $a \geq 1$ there is a number $\lambda_{a}$ such that the limit set of $T_{\lambda}$ has dominated splitting if $\lambda \in (0, \lambda_{a}]$ and does not have dominated splitting for $\lambda \in 
(\lambda_{a},1)$.
 
\subsection{Slap billiard map of an elliptical table}
\begin{figure}
\subfigure[]{
\includegraphics[scale=0.9]{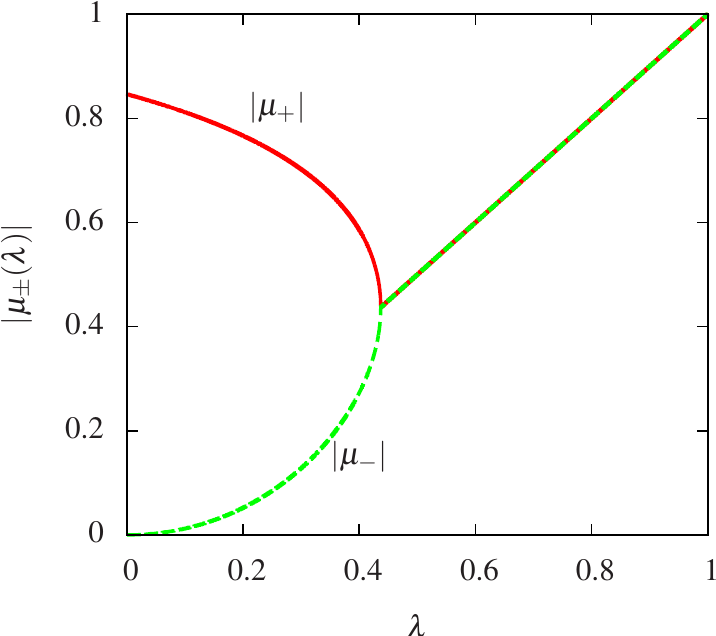}
\label{fig:ellipse-evals}
}
\qquad
\subfigure[]{
\includegraphics[scale=0.9]{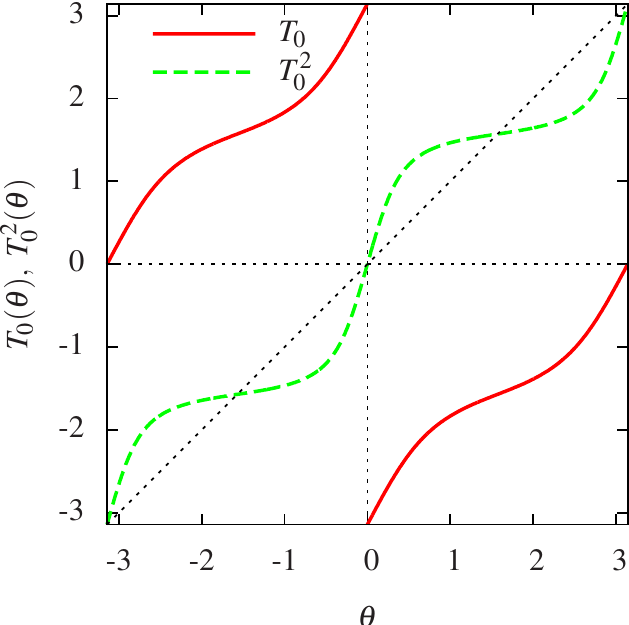}
\label{fig:slap-ellipse}
}
\caption{ (a) Absolute values of the eigenvalues $\mu_{\pm}$ of the stability matrix for the period-$2$ orbit lying along the minor axis in the ellipse with $a=5$.  For $\lambda > \lambda_{-}(5) \simeq 0.44$, the eigenvalues are complex conjugates; for smaller values of $\lambda$, they are real. (b) The slap billiard map $T_0$ and its second iterate $T_0^{2}$ for an
elliptical table with $a=1.2$.
} 
\end{figure}

The slap billiard map $T_0$ for an elliptical table with $a=1.2$ is shown in \figref{fig:slap-ellipse}, together with its second iterate. 
We use the polar angle $\theta$ of the collision point instead of the arc length along the ellipse, since it is simpler to calculate.
In the configuration space (not shown), all trajectories of the slap billiard map are tangent to the evolute of the ellipse, which is a stretched astroid (the Lam\'e curve).

On the graph of $T_0^2$ can be seen the two unstable period-$2$ points lying on the major axis ($\theta = 0$ and $\theta=\pi$), and the two stable period-$2$ orbits on the minor axis. The slap billiard map of any elliptic table is not an expanding map, which agrees with the absence of
chaotic behaviour in the dynamics. 
 
\section{Cardioid}
The results of  \cite{MarkarianPujalsSambarinoPinballBilliards} suggest that pinball billiards may have much more complicated limit sets than just attracting periodic orbits, in particular non-trivial invariant sets with hyperbolic dynamics. From here on we search for more complicated dynamics.

In this section we study pinball billiards in the cardioid, parametrised by the equation
$\rho(\theta) = 1 + \cos(\theta)$
in polar coordinates $(\rho, \theta)$,  with $-\pi \leq \theta < \pi$.
The standard, elastic billiard dynamics in a cardioid has been studied extensively 
\cite{WojtkowskiPrinciplesDesignBilliardsLyapExpsCMP1986, MarkarianNewErgodicBilliardsCardioidNonlin1993, SzaszKPropertyPlanarHypBilliardsCMP1992, BackerDullinSymbDynCardioidBilliardJPA1997, DullinBackerErgodicityLimaconBilliardsNonlin2001},
and in particular has been proved to be hyperbolic, ergodic and mixing.
Robnik \cite{RobnikClassicalBilliardsAnalyticBoundariesJPA1983} studied the existence, stability and bifurcations of periodic orbits in a class of billiards which includes the cardioid, with the focus being the transition from regular dynamics in the circle to stochastic dynamics in the cardioid.
Here we rather study the dynamics in this single table when we vary the contraction parameter $\lambda$ between the one-dimensional map and area-preserving cases.

For the cardioid, the arc length $s(\theta)$ up to an angle $\theta$, measured from the positive $x$-axis, can be evaluated exactly as $s(\theta) = 4 \sin (\theta/2)$ for $0 \leq \theta \leq \pi$. Thus $s=0$ is the rightmost (flattest) point of the cardioid, positive values of $s$ correspond to the top half of the cardioid, and the cusp is at $s=\pm 4$.
This $s$ is used as the abscissa in the phase space plots below.
Some details of the numerical algorithm used to follow the dynamics are given in the appendix.  The curvature of the boundary may be evaluated using the standard expression in polar coordinates:
\begin{equation}
\kappa(\theta) = -\frac{\rho(\theta)^2 + 2 \rho'(\theta)^2 - \rho(\theta) \rho''(\theta)}{\left[\rho(\theta)^2 + \rho'(\theta)^2 \right]^{3/2}}.
\label{eq:curvature-polar}
\end{equation}
The negative sign here corresponds to our convention that focusing boundaries have $\kappa < 0$. The curvature of the cardioid is then $\kappa(\theta) = -3/\sqrt{8(1+\cos \theta)}$.

\subsection{Chaotic attractor}
\begin{figure}
\subfigure[]{
\includegraphics[scale=0.9]{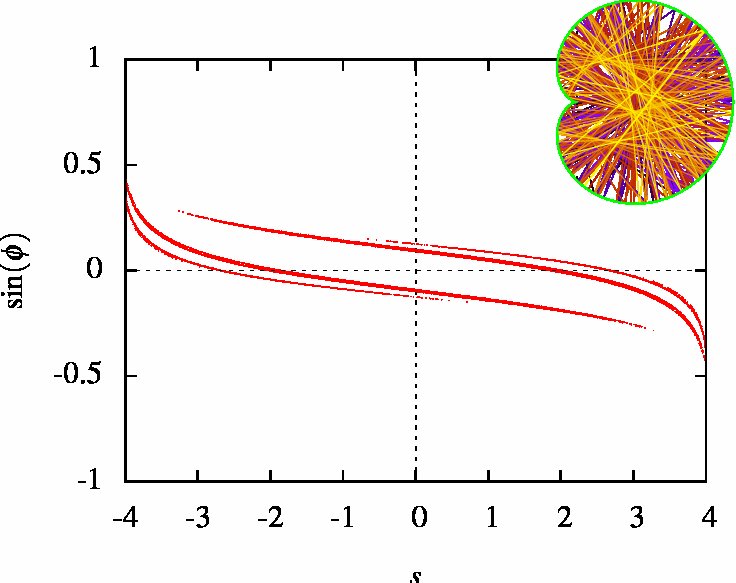}
}
\qquad
\subfigure[]{
\includegraphics[scale=0.9]{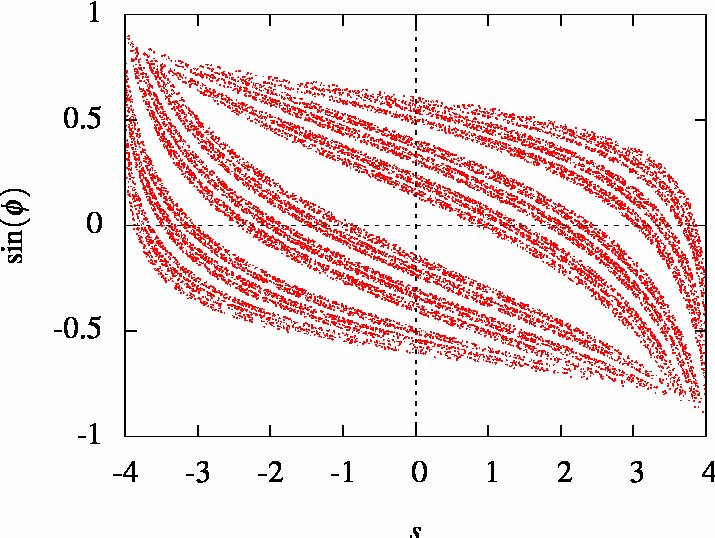}
}
\caption{The chaotic attractor in phase space for the cardioid, for (a) $\lambda = 0.3$ and (b) $\lambda = 0.8$. The inset of (a) shows the attractor in configuration space.
The coordinates are arc length $s$ and $\sin(\phi)$, where $\phi$ is the exit angle at each collision. The cusp of the cardioid is at $s = \pm 4$.}
\label{fig:cardioidAttractor}
\end{figure}

Numerical evidence shows the existence of a non-trivial \emph{chaotic attractor}
for any $\lambda \in [0,1)$, as shown in \figref{fig:cardioidAttractor}. The
attractor  grows in phase space as $\lambda$ increases, until it fills the whole
of phase space uniformly when $\lambda=1$.  Here and in the rest of the paper,
the angular coordinate in phase space is shown as $\sin(\phi)$.

%
  
We remark that points which are arbitrarily close to the singularity $\theta= \pm \pi$ at the cusp of the cardioid (corresponding to $s=\pm 4$) seem to belong to the attractor. However, since the cusp does not belong to the domain of the pinball billiard map, since the normal vector is not uniquely defined there, the attractor does not lie in a compact region of phase space. Thus Theorem 5 of \cite{MarkarianPujalsSambarinoPinballBilliards} cannot be automatically applied. However, the dominated splitting behaviour foreseen by Theorem 4 of that reference is suggested by the following computations.

The trajectory emanating from (almost) any initial condition is found numerically to have two Lyapunov 
exponents $\nu_{\pm}$, with $\nu_+ > 0 > \nu_-$, i.e.\ they are 
uniformly bounded away from zero for any $\lambda \in (0,1]$. The numerical calculation of these Lyapunov exponents is summarised in the appendix, and the results are shown in \figref{fig:cardioidLyapunovExponents}.
Due to the dissipative nature of the dynamics, there is contraction in phase space, which implies that $\nu_+ + \nu_- < 0$.
Each point on the graph, for a fixed value of $\lambda$, is an average over $100$ initial conditions. The Lyapunov exponents evaluated along the different trajectories are found to vary little -- error bars are of the size of the symbols. This, together with visualisations of the attractor for different initial conditions, suggests that there is transitive behaviour on the attractor, i.e.\ that any initial condition on the attractor visits the whole attractor.
These numerical results lead us to
conjecture that the system is in fact \emph{non-uniformly hyperbolic} in the sense of Pesin, i.e.\ almost every point has non-zero Lyapunov exponents \cite[Section~9.2]{RobinsonStabilityDynamicalSystemsBookCRC1998}. However, rigorous proofs of such conjectures are notoriously difficult \cite{RobinsonStabilityDynamicalSystemsBookCRC1998}.

\begin{figure}
\subfigure[]{
\includegraphics[scale=0.9]{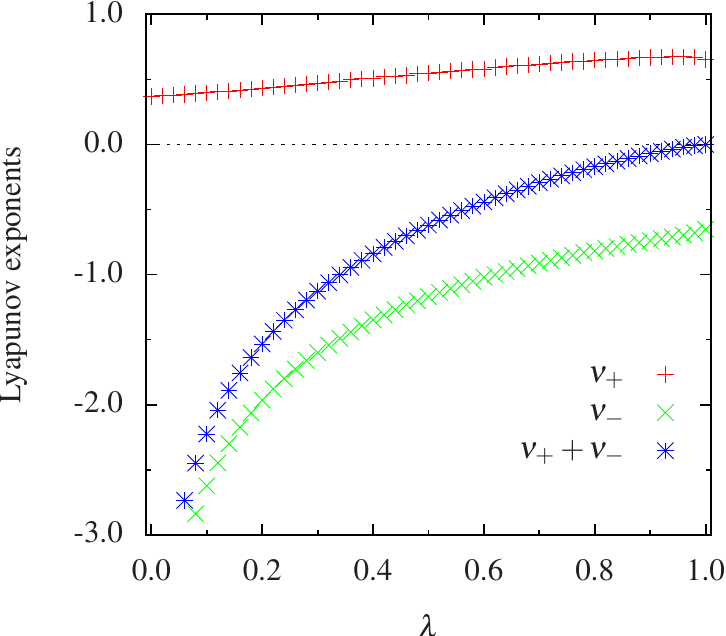}
\label{fig:cardioidLyapunovExponents}
}
\qquad
\subfigure[]{
\includegraphics[scale=0.9]{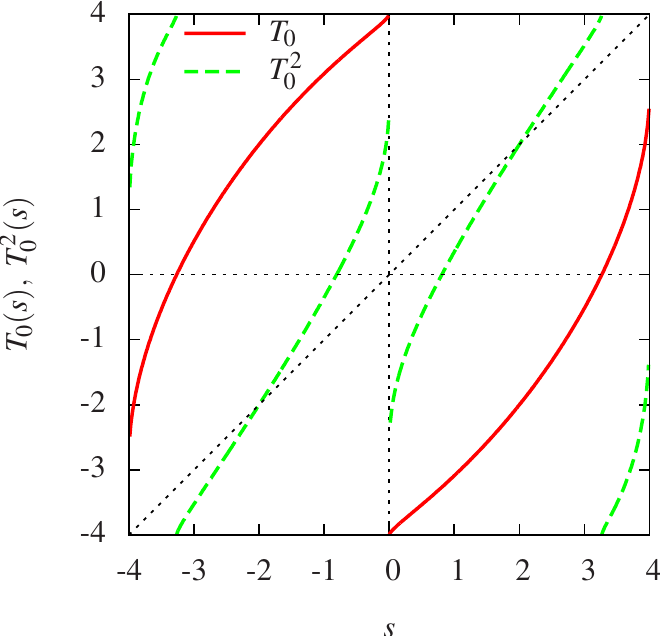}
\label{fig:cardioidSlap}
}
\caption{(a) The two Lyapunov exponents $\nu_{\pm}$ for the cardioid, as a function of $\lambda$. Their sum, which is the logarithm of the area contraction in phase space under the map, is also shown. It converges to $-\infty$ when $\lambda \to 0$, since the map reduces to a one-dimensional map in this case, and to $0$ when $\lambda \to 1$, the elastic case. (b) The slap billiard map $T_0$ for the cardioid and its second iterate. $T_0$ itself has regions where the derivative $T_0'$ is less than $1$ in modulus, near the cusp at $s=0$. 
Nonetheless, its second iterate is expanding, with $\left|T_0'(s)\right|>1$ everywhere, and thus the dynamics of the map are strongly chaotic.
The two unstable fixed points of the second iterate correspond to an unstable vertical period-$2$ orbit joining the top and bottom points of the cardioid at $\theta = \pm \pi/3$, i.e.\ at $s=\pm 2$.
}
\end{figure}
 
\subsection{Slap billiard map of the cardioid}  

In configuration space, all trajectories of $T_{0}$ are tangent to the evolute of the cardioid, which is a mirror image of the cardioid scaled by a factor of $1/3$ (not shown). 
The slap billiard map of the cardioid is interesting in that it has regions,
close to the point $s=0$ of minimal curvature, where the derivative $T_0'$ is
less than one, and other parts where it is greater than one, as shown in
\figref{fig:cardioidSlap}. Nonetheless, as the figure shows, the second iterate
$T_0^2$ of this map is \emph{piecewise expanding}, i.e.\ its derivative is
uniformly bounded away from $1$, with $(T_0^2)'(s) > 1$ everywhere. This
confirms the strong chaotic properties of the dynamics, and suggests a situation
analogous to that in Theorem~12 of 
\cite{MarkarianPujalsSambarinoPinballBilliards} with focusing walls, which
states that under these conditions there should be a finite number of expansive
attractors. In fact, numerically we find that the dynamics of the map seems to
be ergodic.

%
%

\section{Cuspless cardioid}
In the previous two sections, we have studied tables which show completely regular and completely chaotic behaviour, both for elastic and pinball billiard dynamics.
In this section and the following one, we consider two convex tables which exhibit a complicated mixture of phenomena. Here we begin with a convexified version of the cardioid,
obtained by flattening the cusp of the cardioid. To do so, we remove the part of the
cardioid with $\theta \in [-\pi, -\frac{2 \pi}{3}] \cup[\frac{2 \pi}{3}, \pi]$, and replace it by a vertical line at $x=-\frac{1}{4}$ between $y=\pm \frac{1}{4} \sqrt{3}$. This procedure
 gives a convex table which is smooth at all points except for these two joins,
where the curve is only $C^1$ and has discontinuous curvature.   

Again the arc length can be evaluated exactly. We now take $s=0$ as corresponding to the centre of the vertical line at the left of the table, 
and the rightmost point being at $s=\pm C$, where $C \defeq 4 \sin (\pi/3) + \sqrt{3}/4 = 9 \sqrt{3} / 4$. The discontinuities are then at $s = \pm \sqrt{3}/4$.


By unfolding this ``cuspless cardioid'' by reflection in the vertical line, we obtain an equivalent table which satisfies Wojtkowski's defocusing condition \cite{WojtkowskiPrinciplesDesignBilliardsLyapExpsCMP1986, MarkarianBilliardsPesinRegionMeasureOneCMP1988}, and is thus hyperbolic and ergodic in the classical case of elastic collisions \cite{SzaszKPropertyPlanarHypBilliardsCMP1992, MarkarianNewErgodicBilliardsCardioidNonlin1993}.
However, this unfolding process is no longer valid in the case of pinball billiards, since the outgoing angle is modified at each collision with the vertical line.

From the point of view of pinball dynamics, the dynamics always has dominated splitting \cite[Theorem 4]{MarkarianPujalsSambarinoPinballBilliards}. This is verified later, by computing Lyapunov exponents. However, once again the phase space is not compact. Thus the results of \cite{PujalsSambarinoDynamicDominatedSplitting} can be immediately applied only if the attractor is contained in a compact subset of the phase space. This is the case for $\lambda$ close to zero. However, a non-trivial attractor is observed only for larger values of $\lambda$. In fact, it is related to the existence of points of the boundary which are not $C^2$: generic orbits in the chaotic attractor seem to accumulate on these non-smooth points.

\begin{figure}
\subfigure[$\lambda=0.02$]{
\includegraphics{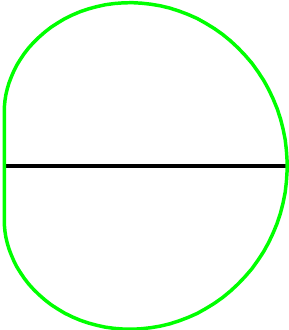}
}
\ 
\subfigure[$\lambda=0.072$]{
\includegraphics{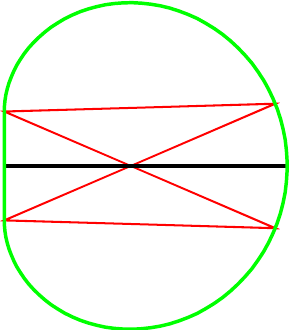}
}
\ 
\subfigure[$\lambda=0.2$]{
\includegraphics{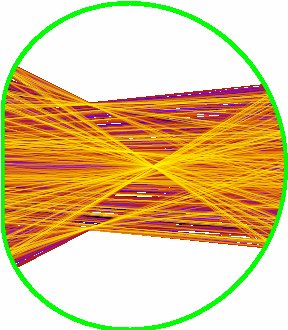}
}
\ 
\subfigure[$\lambda=0.5$]{
\includegraphics{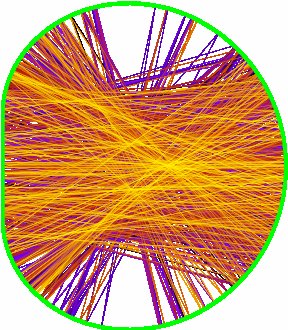}
}
\caption{Evolution of numerically-observed attractors in configuration space with increasing $\lambda$. For $\lambda < \lambda_* \simeq 0.0712$, there is just a period-$2$ attractor.
This periodic attractor coexists with a chaotic attractor for $\lambda \in [\lambda_*, \lambda_c]$, where $\lambda_c \simeq 0.093$. The period-$2$ attractor then becomes unstable, leaving just the chaotic attractor, which expands for increasing $\lambda$.} 
\label{fig:cuspless-attractors-config}
\end{figure}

\begin{figure}
\mbox{
\subfigure[{$\lambda \in [0.1, 0.3]$}]{
\includegraphics[scale=0.9]{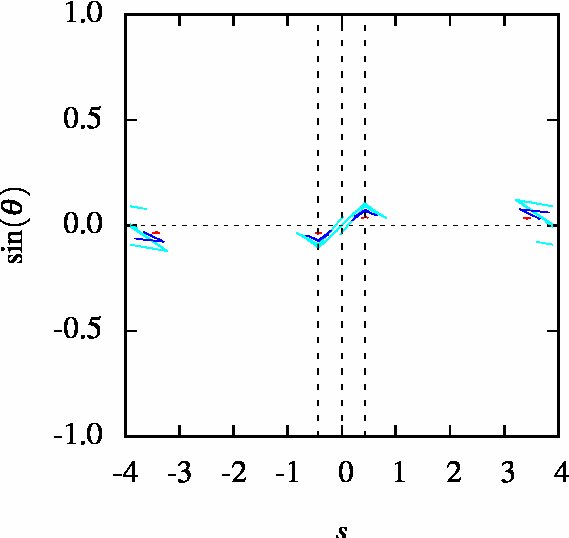}
}
\subfigure[{$\lambda \in [0.4,0.6]$}]{
\includegraphics[scale=0.9]{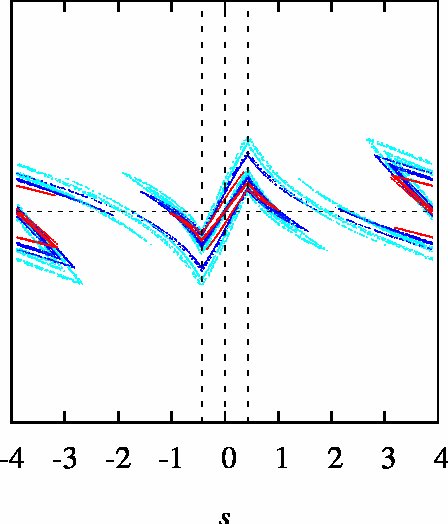}
}
\subfigure[{$\lambda \in [0.7,0.9]$}]{
\includegraphics[scale=0.9]{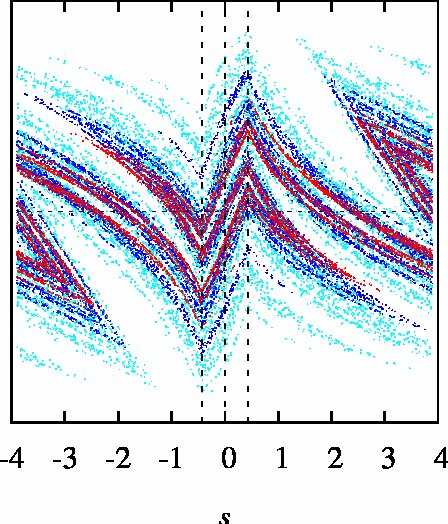}
}
}
\caption{Evolution of the chaotic attractor of the cuspless cardioid in phase space, as a function of $\lambda$. The attractor is depicted for $\lambda$ between $0.1$ and $0.9$ in intervals of $0.1$, with three values of $\lambda$ shown in different colours  in each figure (red, blue and cyan in increasing order of $\lambda$).
Vertical lines mark the centre of the vertical section of the boundary at $s=0$ and the two curvature discontinuities at $s=\pm \sqrt{3}/4$. The attractors are also seen to be non-smooth here; they expand in phase space as $\lambda$ increases. 
}
\label{fig:cardioid-attractors}
\end{figure}

Numerical experiments show that for $\lambda$ close to $0$, the period-$2$ orbit lying in $y=0$ attracts almost all initial conditions, except for the unstable vertical period-$2$ orbit which is the same as in the unmodified cardioid. For $\lambda$ close to $1$, however,  we observe a globally attracting chaotic attractor, as in the unmodified cardioid, although rather different in structure from that case. Careful simulations show that in fact there is a region where the period-$2$ and chaotic attractors \emph{coexist}: some initial conditions are attracted to the periodic orbit, and others to the chaotic attractor. Figures~\ref{fig:cuspless-attractors-config} and \ref{fig:cardioid-attractors} show the evolution of the numerically-determined attractors as $\lambda$ is varied, in configuration space and in phase space, respectively. In the following we shall explain these phenomena.

\subsection{Period-$2$ orbit}
We start by studying the horizontal period-$2$ orbit which runs along the symmetry axis $y=0$, shown in \figref{fig:cuspless-attractors-config}(a). This orbit exists for any $\lambda \in [0,1]$, but its stability depends on the eigenvalues of the stability matrix, i.e.\ the derivative of $T_\lambda ^2$ at the point $x=(0,0)$, which in turn depend on $\lambda$. The stability matrix is obtained by composing the stability matrices at the points $x$ and $T_{\lambda}(x) = (C,0)$, where $C$ is the arc length to the point opposite $0$. We find that
%
\begin{equation}
D_x T_\lambda^2 
= \frac{1}{64}\left(
 \begin {array}{cc}
  -4(27 \lambda + 11)  & 144 (\lambda +1)\\
   33   \lambda (\lambda +1)&   - 4 \lambda (11 \lambda + 27)
\end{array}
\right). 
\end{equation}
The eigenvalues $\mu$ of this matrix are
\begin{equation}
 \mu_{\pm}(\lambda) = \frac{1}{32} \left(-11-54 \lambda -11 \lambda ^2 \pm \sqrt{11} \
(1+\lambda ) \sqrt{11+86 \lambda +11 \lambda ^2}\right).
\end{equation}
Both eigenvalues are real and satisfy $\mu_\pm(\lambda) \leq 0$, for any $\lambda \in [0,1]$. Moreover, $\mu_\pm(\lambda) > -1$ for any $\lambda \in (0, \lambda_{c})$, where $\lambda_{c} \defeq \frac{1}{5}(27 - 8 \sqrt{11}) \simeq 0.0934$ is the value such that that $\mu_-(\lambda_{c}) = -1$, i.e.\ for which this eigenvalue becomes unstable. 
 The stable orbit is visualised later for the slap map $T_{0}$ (see \figref{fig:cusplessSlap}).

When an eigenvalue crosses $-1$, a period-doubling bifurcation is expected \cite[Chapter~7]{RobinsonStabilityDynamicalSystemsBookCRC1998}. 
However, in this case no \emph{stable} period-$4$ orbit emerges, which would be observed as an attractor. We expect that an \emph{unstable} period-$4$ orbit is created in this bifurcation; presumably it could be rendered stable if the billiard table were modified to have smaller curvature around $s=0$.


\subsection{Birth of chaotic attractor}

In numerical experiments, a globally attracting chaotic attractor is observed for $\lambda \gtrsim 0.1$ where the period-$2$ orbit is no longer stable. However, this chaotic attractor is also observed to persist for  $\lambda > \lambda_{*} \simeq 0.0712$ in \emph{coexistence} with the attracting period-$2$ orbit discussed above.

\begin{figure}
\subfigure[]
{
\includegraphics[scale=0.9]{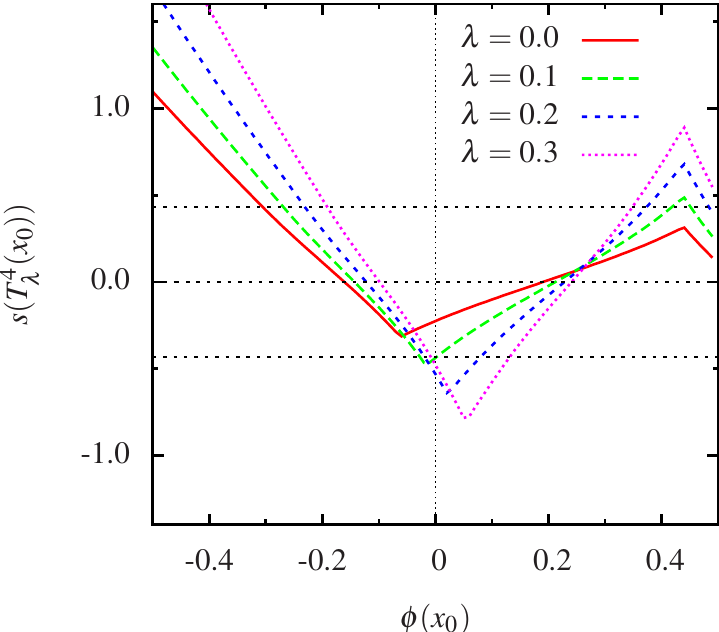}
\label{fig:cuspless-period-4}
}
\quad
\subfigure[]
{
\includegraphics[scale=0.9]{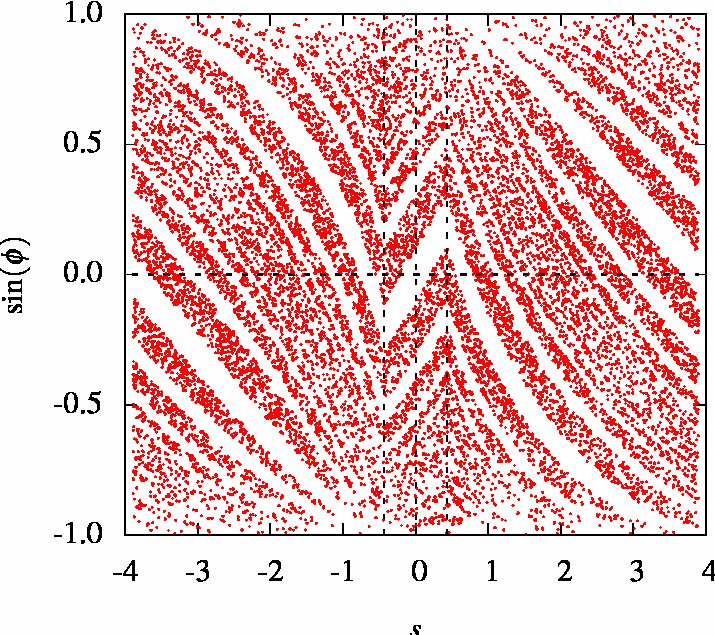}
\label{fig:cardioid-basins}
}
\caption{(a)  Arc length $s_4$ of the fourth collision point, starting from initial conditions $x_0 = (s_0, \phi_0)$ at the upper discontinuity $s_0 \defeq -\sqrt{3}/4$ in the cuspless cardioid, as a function of the initial angle $\phi_0$. The dotted horizontal lines show the two discontinuities $s=\pm \sqrt{3}/4$ and $s=0$ (the centre of the vertical line). 
(b) Basin of attraction of the chaotic attractor for the cuspless cardioid with $\lambda=0.09$, shown in phase space with coordinates $s$ and $\sin(\phi)$, where $\phi$ is the exit 
 angle. The points which are not attracted to the chaotic attractor (shown in blank) are attracted to the stable period-$2$ orbit. 
The basins appear to be intermingled in some parts of phase space.
} 
\end{figure}

How and why does this chaotic attractor arise? For $\lambda$ just above $\lambda^*$ it is observed to be very thin, and is concentrated around what appears to be a periodic orbit of period~$4$  joining the two points of discontinuity in the curvature (figure~\ref{fig:cuspless-attractors-config}(b)).
We search for such a period-$4$ orbit as follows. Starting exactly at the upper discontinuity, we shoot trajectories with different initial angles $\phi$, and follow them for four collisions. Plotting the arc length $s$ of the collision point at the fourth collision gives \figref{fig:cuspless-period-4}, for angles close to $\phi=0$, which are those for which the attractor appears. 

We see that for $\lambda < \lambda^*$, the fourth collision occurs for values of arc length \emph{greater} than $s^*$, that is, the collisions occur on the vertical line joining the two discontinuities. The trajectory then converges towards the attracting period-$2$ orbit.  For $\lambda>\lambda^*$, the trajectory returns to locations with $s < s^*$, that is it collides on the curved surface just above the discontinuity. Nonetheless, due to the focusing nature of the boundary, the trajectory is then reinjected into the same region, and in this case we observe that the subsequent dynamics converges to a chaotic attractor for any $\lambda > \lambda^*$.  

For $\lambda = \lambda^*$, the trajectory returns exactly to the same discontinuity on the fourth collision, and there is indeed a period-$4$ orbit joining the two discontinuities exactly at this value of $\lambda$. 
By searching for the pair $(\lambda, \phi)$ for which the curve of \figref{fig:cuspless-period-4} just touches $s^*$, we estimate the critical values as $\lambda^* \simeq 0.0712173$ and $\phi^* \simeq -0.0289796$ (see figure \ref{fig:cuspless-period-4}). 
Simulations confirm the existence of the orbit and its instability.

Since the period-$2$ and chaotic attractors coexist in the interval $\lambda \in [\lambda^*, \lambda_c]$, it is of interest to examine their respective \emph{basins of attraction} in phase space, i.e.\ the sets of initial conditions in phase space which converge to each attractor.
These basins are shown in \figref{fig:cardioid-basins} for the case $\lambda = 0.09$. In certain regions the basins appear to be intermingled arbitrarily closely.

%

\subsection{Chaotic attractor}
As mentioned above, a chaotic attractor is found for all $\lambda > \lambda^*$, and for all $\lambda > \lambda_c$ this is the only attractor.
Numerical results again indicate that 
there is one positive and one negative Lyapunov exponent which are both constant almost everywhere, so that the chaotic attractor again appears to be non-uniformly hyperbolic.

\begin{figure}
\subfigure[]{
\includegraphics[scale=0.9]{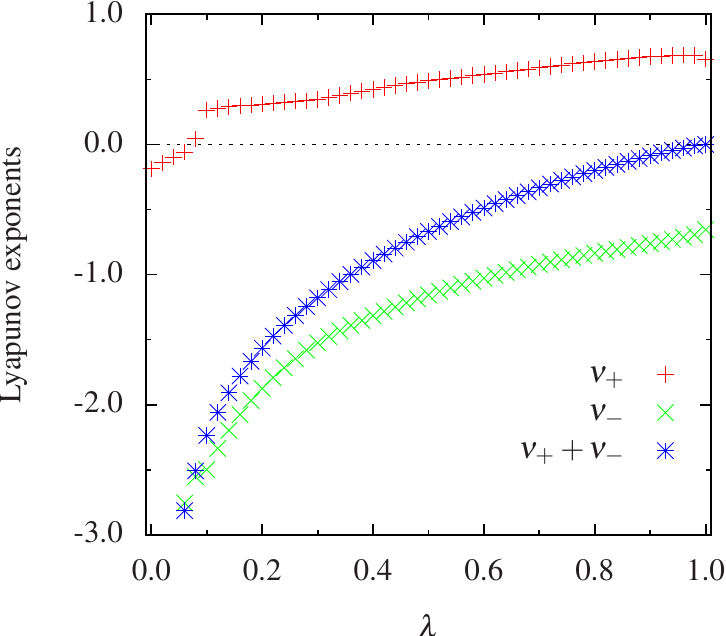}
\label{fig:cusplessLyapunovExponents}
}
\qquad
\subfigure[]{
\includegraphics[scale=0.9]{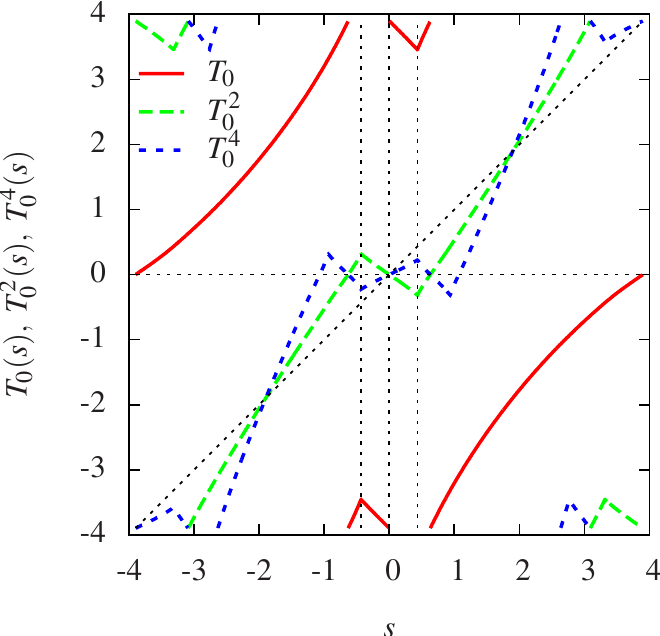}
\label{fig:cusplessSlap}
}
\caption{(a) Lyapunov exponents for the cuspless cardioid, averaged over $100$ independent trajectories. 
Note that $\nu_+$ crosses $0$ when the chaotic attractor
appears, and has a jump where the period-$2$ orbit becomes unstable, since then there are no longer attracting periodic orbits with $\nu_+ < 0$. (b) Slap billiard map, together with its second and fourth iterates, for the cuspless cardioid.}
\end{figure}

As $\lambda$ increases, the attractor expands in phase space, as shown 
in \figref{fig:cardioid-attractors}, until it fills the whole of phase space when $\lambda=1$.  In particular, it visits a larger and larger part of the boundary.
Note that the ``corners'' which are visible in plots of the attractors correspond to the positions of the discontinuities. Thus the non-differentiability of the map gives rise to a
non-smooth shape of the attractor. 
The form of the attractor is reminiscent of the H\'enon attractor \cite{HenonTwoDimensionalMappingStrangeAttractor1976} and its linear counterpart, the Lozi map \cite{LoziAttractor1978}. This gives an example where the conclusion of \cite{MarkarianPujalsSambarinoPinballBilliards} on the existence of complicated dynamics
holds, but with singularities which are not accounted for in that reference.

\subsection{Slap billiard map for the cuspless cardioid}
Figure \ref{fig:cusplessSlap} shows the slap billiard map and its second and fourth iterates for the cuspless billiard.
The graph of the 4th iterate confirms that there are no period-$4$ orbits. 


\section{Smooth deformations of the circle: three-pointed egg}

The last example of focusing pinball billiards that we consider is a family of smooth deformations of the circle, given by the equation
\begin{equation}
 \rho(\theta) = 1 + \alpha \cos(3 \theta)
\end{equation}
in polar coordinates $(\rho, \theta)$ with $\theta \in [-\pi, \pi)$;
they are related to a class of
billiards studied by Robnik in the case of elastic collisions \cite{RobnikClassicalBilliardsAnalyticBoundariesJPA1983}.   
The shape parameter $\alpha$ determines the exact shape of the table: $\alpha=0$ gives a circular table, while small values of $\alpha$ give smooth perturbations of the circle with $3$-fold rotational symmetry, a shape which we call a ``three-pointed egg'', shown in \figref{fig:period3}.
However, when $\alpha > 1/10$ the table becomes non-convex and non-focusing, since ``waist'' regions develop;  we thus consider  only $\alpha \in [0,1/10]$ in this paper. An important feature of this class of tables is that the boundary is analytic, and hence contains no singular points.
The domain of the pinball billiard map is thus compact.

As in the ellipse, the arc length is difficult to calculate, and we again replace it by the angle $\theta$ of collision with respect to the positive $x$-axis.
The curvature may be evaluated using \eqref{eq:curvature-polar}.
The minimum curvature occurs for $\check{\theta}_{0} \defeq \pi/3$ and its symmetric
counterparts $\check{\theta}_{1} \defeq \pi$ and $\check{\theta}_{2} \defeq -\pi/3$; it is
\begin{equation}
\kappa(\check{\theta}_{i}) = -\frac{1 - 10 \alpha}{(1-\alpha)^2}  \qquad \mathrm{ for\ } i = 0,1,2.
\end{equation}
We denote by $\hat{\theta}_{i} \defeq \check{\theta}_i + \pi\ (\mathrm{mod\ } 2\pi)$ the boundary positions diametrically opposite the $\check{\theta}_{i}$; these are the points of maximal curvature, with $\kappa(\hat{\theta}_i) = -\frac{1+10 \alpha}{(1+\alpha)^2}$.

\subsection{Hamiltonian case, $\lambda=1$}

\begin{figure}
\subfigure[$\alpha=0.02$]
{
\includegraphics[scale=0.9]{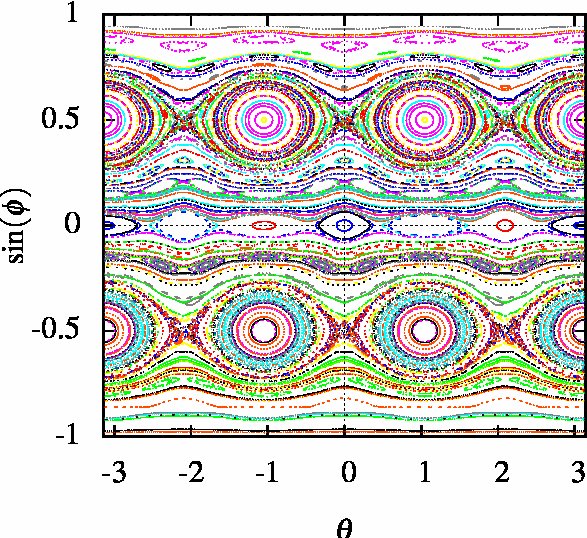}
}
\subfigure[$\alpha=0.05$]
{
\includegraphics[scale=0.9]{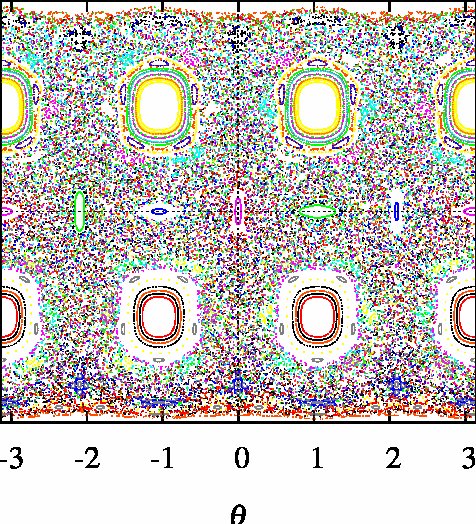}
}\subfigure[$\alpha=0.08$]
{
\includegraphics[scale=0.9]{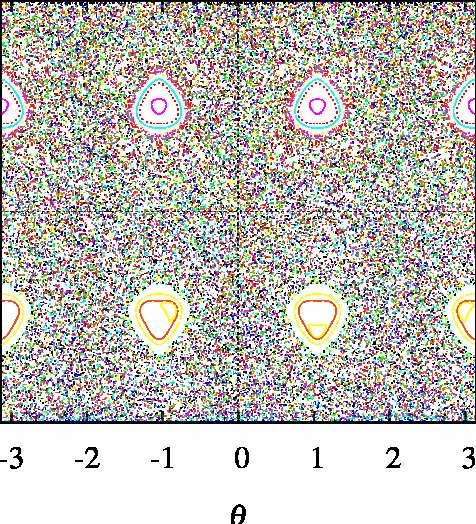}
}
\caption{Phase space of the three-pointed egg in the Hamiltonian case $\lambda = 1$
(elastic collisions), for three values of the shape parameter $\alpha$ as indicated.
Different colours indicate 
trajectories from different initial conditions.
The phase spaces are of mixed type,
consisting of elliptic islands surrounding periodic orbits, together with a chaotic sea in (b) and (c).
The elliptic period-$3$ orbits can be clearly seen in all three cases along the lines $\sin(\phi) = \pm 0.5$, and the hyperbolic ones can be seen in (a) on the same line. In (a) and (b), the elliptic period-$2$ orbits can also be observed on the line $\phi =0$; these have become unstable for the value of $\alpha$ in (c).
}
\label{fig:islands}
\end{figure}

It is informative to first study the Hamiltonian case, when $\lambda=1$ and the 
collisions are elastic. In this case, for any $\alpha \in [0, 1/10]$ the system exhibits a generic mixed-type phase space, consisting of elliptic islands and a chaotic sea, as shown in \figref{fig:islands} for three values of the parameter $\alpha$ in the allowed range. 

For $\alpha = 0$, the dynamics is completely integrable, and the phase space is foliated by horizontal invariant curves. As $\alpha$ increases, alternating elliptic and hyperbolic periodic points appear \cite[Section~3.2]{LichtenbergLiebermanRegularChaoticDynamicsBook2ndEd1992}. Due to the $3$-fold symmetry of the tables, there is a bias towards period-$3$ orbits: 
the figures exhibit elliptic points of period $3$, with their associated islands, and for small enough $\alpha$ also islands around elliptic period-$2$ points. For larger values of $\alpha$, there is a transition towards Hamiltonian chaos with a chaotic sea \cite[Chapter~4]{LichtenbergLiebermanRegularChaoticDynamicsBook2ndEd1992}, although here full chaos is never reached.

\subsection{Slap billiard map}

\begin{figure}
\subfigure[]
{
\includegraphics[scale=0.9]{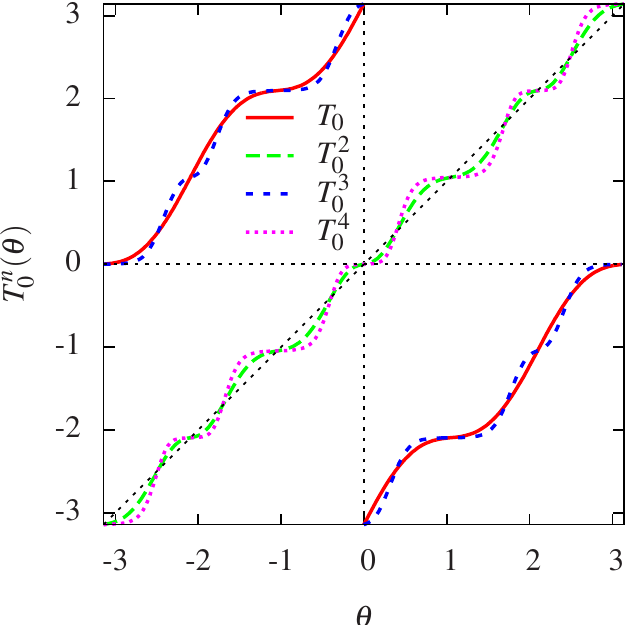}
}
\qquad
\subfigure[]
{
\includegraphics[scale=0.9]{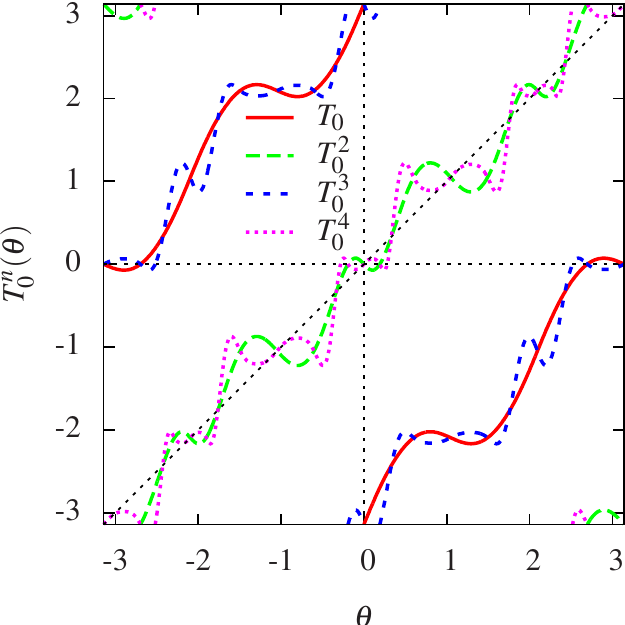}
}
\caption{Iterates up to the fourth of the slap billiard map ($\lambda = 0$) for the three-pointed egg, for (a) $\alpha = 0.05$ and (b) $\alpha = 0.08$.
In (b), critical points have developed, and the previously stable period-$2$ orbits have become unstable in a period-doubling bifurcation, giving rise to stable period-$4$ orbits nearby.
There are no period-$3$ orbits.
} 
\label{fig:polar-slap}
\end{figure}

To understand the dynamics of this class of billiards, we next study the slap billiard map $T_0$, which gives information on the periodic orbits which occur when $\lambda = 0$. $T_0$ is shown in \figref{fig:polar-slap}, as a function of the angle $\theta$ around the origin, together with its iterates up to the fourth. When the shape parameter $\alpha$ tends to zero, $T_0$ converges to the translation $x + \pi$, while for $\alpha$ larger than a critical value, three critical points arise, i.e.\ $\theta$ for which $T_0'(\theta)=0$. This implies that there are pairs of boundary points whose normals have the same intersection point on the opposite side of the table.

The graph of $T_0^3$ shows that there are no period-$3$ orbits when $\lambda=0$. Thus the period-$3$ orbits which are found for $\lambda = 1$ are created at some non-trivial value of $\lambda \in (0,1)$.  For $\alpha = 0.05$ there are stable period-$2$ orbits joining the vertices of the billiard with their opposite flat regions. These have become unstable by $\alpha = 0.08$, and so are no longer visible.
There are also additional period-$2$ orbits, which are always unstable -- these join parts of the billiard at which the tangents are parallel. 
The following sections study all of these classes of orbits in more detail.

\subsection{Period-$2$ orbits}
\label{sec:polar-period-2}
Each vertex $\hat{\theta}_i$ of the three-pointed egg is joined to its diametrically-opposite point of minimum curvature $\check{\theta}_i$ by a periodic orbit of period $2$, with outgoing angle zero; these orbits are depicted in \figref{fig:period3}. They exist for any $\alpha$, but can become unstable. Their stability matrix, which depends on $\alpha$ and $\lambda$, is given by
\begin{equation}
\left(
\begin {array}{cc} 
\frac {1+\alpha^4-2 \alpha^2  \left( 163+162 \lambda \right) }{ \left(\alpha^2 - 1 \right) ^2} 
&
\frac { 2( \alpha^2 + 18 \alpha - 1 ) ( 1+\lambda ) }{ ( \alpha - 1 ) ^2} 
\\ 
-\frac {162 \alpha^2 ( \alpha^2-18 \alpha-1  ) \lambda(\lambda + 1 ) } 
{ (\alpha - 1 ) ^2 (\alpha + 1 ) ^4}
&
\frac {\lambda\,  ( \lambda^2+\alpha^4 \lambda^2 - 2 \alpha^2 \lambda  ( 162+163 \lambda  )  ) }
{ \left( \alpha^2 - 1 \right) ^2}
\end {array} 
\right).
\end{equation}

Let $\mu_{\pm}(\alpha, \lambda)$ denote the eigenvalues of the stability matrix.
For $(\alpha, \lambda) \in [0,0.1] \times [0,1]$ we have
$\mu_{\pm}(\alpha, \lambda)\neq 1$. However, 
defining
\begin{equation}
   \tilde \alpha(\lambda) \defeq \left[\frac{82+162 \lambda +82 \lambda ^2-9 (1+\lambda ) 
\sqrt{83+162 \lambda +83 \lambda ^2}}{1+\lambda ^2} \right]^{1/2},
\end{equation}
we have $\mu_{-}(\tilde \alpha(\lambda), \lambda) = -1$, i.e.\ one of the eigenvalues crosses $-1$, giving rise to a period-doubling bifurcation: the period-$2$ orbits become unstable for $\alpha > \tilde \alpha(\lambda)$, and a stable period-$4$ orbit is created nearby in phase space. If $\alpha > \tilde \alpha(0) \simeq 0.0781$, then  the period-$2$ orbit is unstable for all $\lambda$. The curve $\tilde \alpha(\lambda)$ is shown in \figref{fig:polar-phase-diag}.
 
On the other hand, for $\alpha <  \tilde{\alpha}(\lambda)$,
both eigenvalues have modulus less than one, so the period-$2$
orbits are stable. Moreover, for pairs $(\alpha, \lambda)$ such that $\alpha <
\tilde \alpha(1)\simeq 0.0554$ and
\begin{equation}
 \alpha > \hat \alpha(\lambda) \defeq \frac{9(1+\lambda) - \sqrt{82 + 160 \lambda + 82 \lambda^2}}{\lambda-1},
\end{equation}
the eigenvalues are complex conjugates. 
 


\subsection{Period-$3$ orbits}

For elastic collisions ($\lambda=1$), there are two elliptic period-$3$ orbits of non-zero angle which join the flat parts of the boundary, $\{\check{\theta}_{i}\}$, as can be seen in \figref{fig:islands}.
These orbits come in pairs, with the same points on the boundary but opposite angles, and by symmetry they form equilateral triangles in configuration space. It can be verified that they are elliptic (i.e.\ the eigenvalues of the stability matrix lie on the unit circle) for any $\alpha$.

When $\lambda$ is decreased away from $1$,  these period-$3$ orbits persist, but they become \emph{attracting}, a phenomenon which is well known for weak dissipative perturbations of Hamiltonian systems \cite{LiebermanTsangTransientChaosDissipativelyPerturbedHamiltonianPRL1985, FeudelGrebogiMap100CoexistingPeriodicAttractorsPRE1996,GallasMultistabilityKickedRotorIJBC2008}.
Higher-order periodic attractors are also expected for very weak dissipation \cite{FeudelGrebogiMap100CoexistingPeriodicAttractorsPRE1996}, and these are indeed also found in our system, although we have not investigated this point in detail.

\begin{figure}
\subfigure[]
{
\includegraphics[scale=0.9]{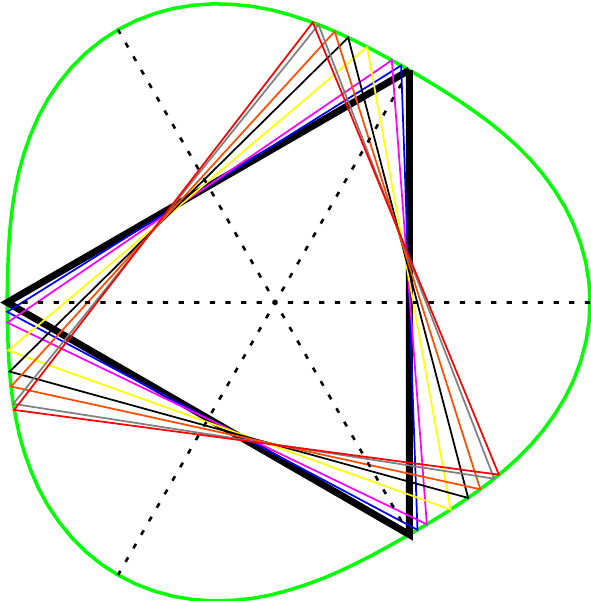}
\label{fig:period3}
}
\qquad
\subfigure[]
{
\includegraphics[scale=0.9]{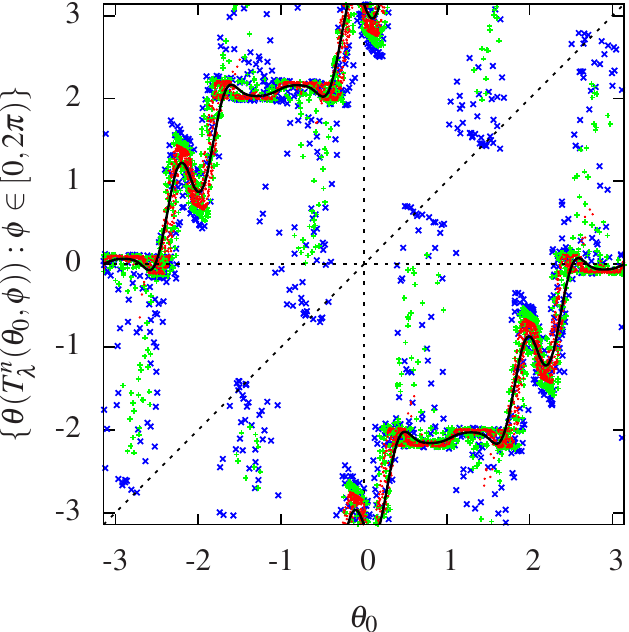}
\label{fig:polar-nonslap}
}
\caption{(a) Change in position of one of the attracting period-$3$ orbits as
$\lambda$ is varied, with the shape parameter $\alpha=0.08$. For
$\lambda=1$ (thick black line) the orbit is elliptic; for $\lambda <
1$ it is attracting. The values of $\lambda$ shown are, in an
anti-clockwise direction, $\lambda = 1.0, 0.9, 0.8, 0.7, 0.6, 0.5,
0.45, 0.41, 0.40$. The latter value is close to the
numerically-determined limit of existence of the period-$3$ orbits,
which is $\lambda \simeq 0.39$. For this value of $\alpha$, the period-$2$
orbits are unstable; they are shown as dotted lines. 
(b) Third iterate of the perturbed slap map for the three-pointed egg with $\alpha=0.08$, for $\lambda=0.0$ (slap map; solid back line), $\lambda = 0.2$ (red dots), $\lambda = 0.3$ (green $+$), and $\lambda = 0.4$ (blue $\times$).
The fact that there is no intersection with the diagonal for $\lambda < 0.4$ shows that there are no period-$3$ orbits in those cases. 
} 
\end{figure}

The symmetry of the table implies that the orbits remain equilateral triangles, but they must rotate away from their initial position in order to accommodate the altered reflection angles.  \bfigref{fig:period3} shows the change in position of one of these period-$3$ orbits as $\lambda$ varies in the case $\alpha = 0.08$.

Below a value $\lambda = \lambda_c(\alpha)$, which depends on $\alpha$, these orbits cease to exist. To show this we use the following argument.
Recall that the graph of $T_{0}^3$ does not intersect the graph of the identity, proving that there are no periodic orbits of period~3 for the slap map. This argument extends to $T_{\lambda}^3$ for small values of $\lambda>0$, as follows. 
Consider the projection $\pi(\theta,\phi) = \theta$ to the space of positions only, and study the set $\{ ( \theta_0, \pi \circ T_{\lambda}^3(\theta_0,\phi): (\theta_0,\phi) \in M \}$. This consists of all pairs $(\theta_0, \theta_3)$, where $\theta_0$ is an initial position and $\theta_3$ is the position of the third iterate of the initial point, chosen with \emph{any} initial angle. If the graph of this set
 does not intersect the graph of the identity map $\{(\theta,\theta)\}$, then we can conclude that there are \emph{no} periodic orbits of period~3. 
Of course, if there is such an intersection, then we can conclude nothing -- in particular, this does \emph{not} imply that there are such period-$3$ orbits.  

\bfigref{fig:polar-nonslap} shows this map for $\alpha = 0.08$ and small values of $\lambda$, from which we indeed conclude that $T_\lambda$ has no periodic orbits of period $3$ below some $\lambda_c$, with $\lambda_c \simeq 0.39$ in this case.  
For $\lambda>\lambda_c$, there are intersections, but we reiterate that no firm conclusions can be drawn from this; nonetheless, attracting period-$3$ orbits are found numerically for all $\lambda \geq \lambda_c$.

These results suggests that the period-$3$ orbits disappear in a saddle--node bifurcation at $\lambda=\lambda_c$, by annihilating with another class of period-$3$ orbits. This latter class of orbits joins the vertices $\hat{\theta}_i$ of the three-pointed egg when $\lambda=1$, and are always unstable. They also rotate upon changing $\lambda$. We conjecture that these orbits collide with the stable period-$3$ orbits at $\lambda = \lambda_c$. Further evidence for this is obtained by calculating the maximal Lyapunov exponent of the stable period-$3$ orbits as $\lambda \to \lambda_c$. The Lyapunov exponent converges to $0$, suggesting an eigenvalue $1$ of the stability matrix, which indeed corresponds to a saddle--node bifurcation \cite[Chapter~7]{RobinsonStabilityDynamicalSystemsBookCRC1998}.

\subsection{Global picture of parameter space}

As we vary the parameters $(\alpha, \lambda)$ in the two-dimensional parameter space $[0,0.1] \times [0,1]$, we numerically observe a variety of different types and combinations of attractors. Roughly, for $\alpha$ close to $0$ there are only periodic attractors, while for $\alpha$ close to $0.1$ there is a range of $\lambda$ for which chaotic attractors are observed to coexist with periodic ones, in particular with the period-$3$ attractors described above. 

\begin{figure}
\includegraphics{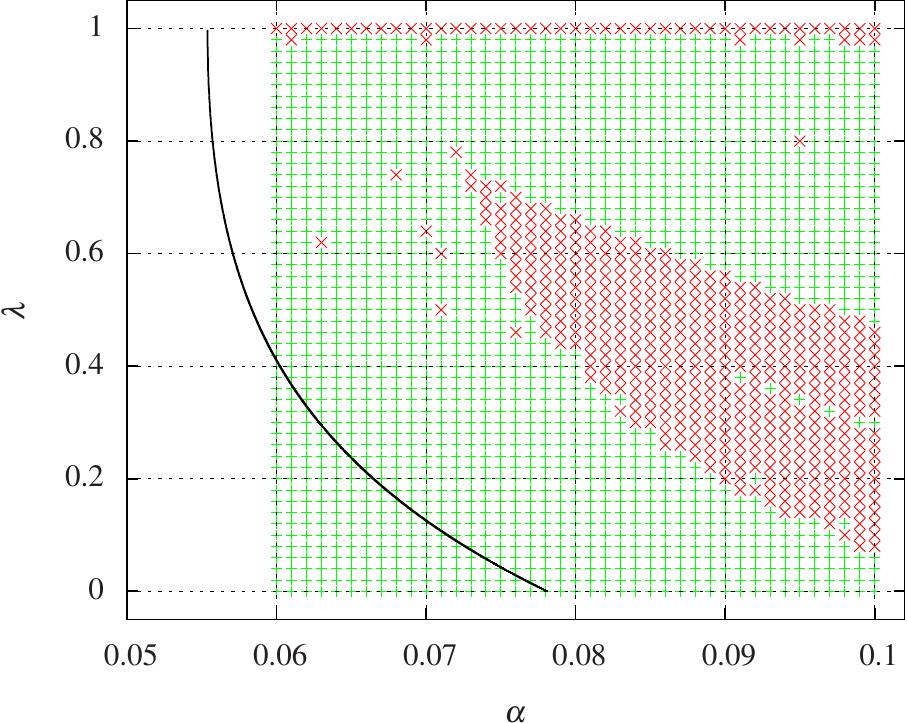}
\caption{Parameter space of the three-pointed egg, showing the regions of parameters  where some (red $\times$) or no (green $+$) chaotic behaviour is found. This was calculated by taking $100$ random initial conditions for each $(\alpha, \lambda)$ pair and calculating the corresponding Lyapunov exponents. If any of the trajectories has a positive Lyapunov exponent, then $(\alpha, \lambda)$ is assigned to the region with chaotic behaviour. 
No chaotic attractors were found for $\alpha < 0.06$, except for the elastic case $\lambda=1$, where there is always some chaotic behaviour for any $\alpha>0$.
Any attractors whose basins occupy only a very small fraction of phase space are
not seen in this calculation, and thus do not appear in the figure.  The black
curve represents the function $\tilde{\alpha}(\lambda)$ calculated in
section~\ref{sec:polar-period-2}, which 
gives the limit of stability of the period-$2$ orbits: above the curve, the
period-$2$ orbits are unstable; below it they are stable.
}
\label{fig:polar-phase-diag}
\end{figure}

In order to begin to understand this phenomenology, we first give in \figref{fig:polar-phase-diag} a global picture of parameter space,  which depicts only the most gross level of information: whether or not chaotic behaviour is found for a given $(\alpha, \lambda)$ pair, that is whether there are any initial conditions which lie in the basin of attraction of some chaotic attractor, characterised by having a positive Lyapunov exponent. We can conclude from the figure, coupled with the analysis of the stability of period-$2$ orbits, that increasing $\alpha$ tends to destabilise the dynamics, as could be expected. 

We remark that for $\lambda=1$ and any $\alpha>0$ there are always some chaotic orbits, namely those which wander in the chaotic sea exhibited in \figref{fig:islands} (or those which stay confined to a stochastic layer for smaller values of $\alpha$, which is the reason why $\lambda=1$ is indicated as showing chaotic behaviour in \figref{fig:polar-phase-diag}. The red crosses at $\lambda=1$ should thus extend over the whole range of allowed $\alpha$.

Since the billiard table, and hence the pinball billiard map $T_\lambda$, are smooth, we can invoke the implicit function theorem to show that periodic orbits, both attracting and hyperbolic, should usually be expected to persist under small enough variations both of $\alpha$ and of $\lambda$. We are thus especially interested in the situations where this does not hold, that is when bifurcations occur. The persistence of chaotic attractors under perturbation is less clear, since these can undergo different types of \emph{crisis}, where they suddenly change their extent or indeed disappear completely due to the interaction with other features of the dynamics \cite[Chapter~8]{OttChaosDynamicalSystemsBook1993}.

\subsection{Period-doubling cascade}

The main feature visible in \figref{fig:polar-phase-diag} is the ``tongue'' region for large $\alpha$ where chaotic attractors occur. To understand this, we consider a cut across parameter space by fixing $\alpha=0.08$ and studying the evolution of the numerically-observed attractors as $\lambda$ is varied. This gives representative results for other values of $\alpha$ in this region.

As discussed above, for this value of $\alpha$, there are attracting period-$4$ orbits for $\lambda=0$, and numerically these are found also for small $\lambda$. For $\lambda > 0.39$, there are attracting period-$3$ orbits, as shown above. But \figref{fig:polar-phase-diag} shows that there is also a region of values of $\lambda$, approximately the interval $[0.42, 0.66]$, which overlaps the region of existence of attracting period-$3$ orbits, where chaotic attractors are also found.  It is thus of interest to investigate the route by which the chaotic behaviour is first created and then destroyed.

\begin{figure}
\mbox{
\subfigure[$\lambda=0.1$]
{
\includegraphics{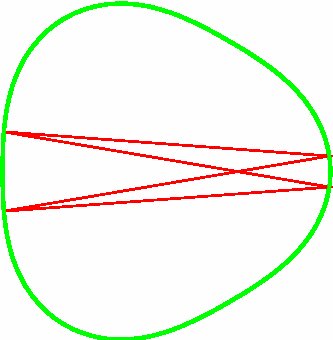}
}
\subfigure[$\lambda=0.39$]
{
\includegraphics{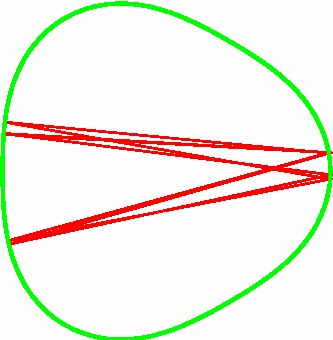}
}
\subfigure[$\lambda=0.43$]
{
\includegraphics{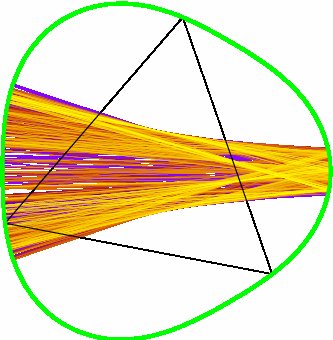}
}
\subfigure[$\lambda=0.45$]
{
\includegraphics{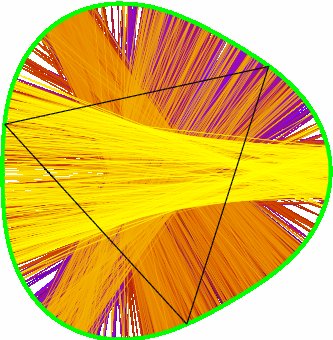}
}
}
\caption{Attractors in configuration space for the three-pointed egg with $\alpha=0.08$: (a) period-$4$ orbit which persists from the stable period-$4$ orbit of the slap map ($\lambda=0$); (b) period-$8$, after undergoing a single period-doubling bifurcation; (c): localised chaotic attractor, after the accumulation of period doublings; (d): delocalised chaotic attractor, after the merging crisis. Here, the trajectories tend to remain for a long time in each part of the attractor that was previously localised, before jumping to a different part, as shown by the colours in the figure. In each of (c) and (d), a coexisting period-$3$ orbit is also shown in black. Note that for this value of $\alpha$, the horizontal period-$2$ orbit is unstable even for $\lambda=0$, and thus does not attract generic initial conditions. } \label{fig:polar_attractors_config}
\end{figure}

\begin{figure}
\subfigure{
\includegraphics[scale=0.9]{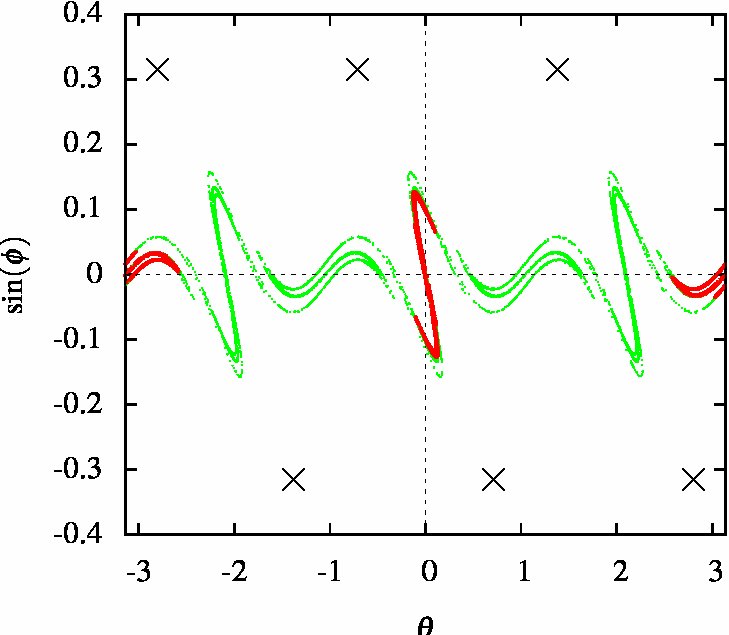}
}
\quad
\subfigure{
\includegraphics[scale=0.9]{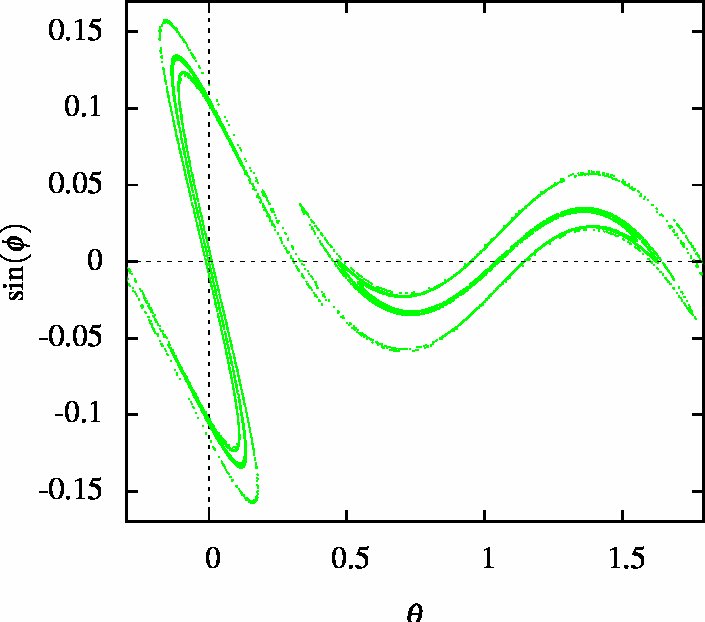}
}
\caption{(a) Attractors in phase space for the three-pointed egg with $\alpha=0.08$ and $\lambda=0.43$ (red), $\lambda=0.45$ (green). The coexisting period-$3$ attractors are also shown (black $\times$). (b) Enlargement of the attractor for $\lambda=0.45$.} 
\label{fig:polar-attractors-phase}
\end{figure}

Numerically we find that this chaotic behaviour arises by the well-known period-doubling route \cite[Chapter~8]{OttChaosDynamicalSystemsBook1993}, where a cascade of an infinite number of successive period doublings at increasing values of $\lambda$ accumulates at a certain value $\lambda_\infty$. For $\lambda>\lambda_\infty$, there is a chaotic attractor with a positive Lyapunov exponent, except in windows where periodic behaviour is observed, similar to the periodic windows observed in the logistic map  \cite[Chapter~2]{OttChaosDynamicalSystemsBook1993}. The development of the observed attractors is shown in configuration space in \figref{fig:polar_attractors_config} and in phase space in \figref{fig:polar-attractors-phase}.

For small values of $\lambda$, there are only stable period-$4$ orbits.
At a critical value of $\lambda$, these become unstable via a
period-doubling bifurcation, giving stable period-$8$ orbits. As $\lambda$ increases, we observe a repeated doubling of
the period of the attracting orbits, giving orbits with period $2^n$ at $\lambda_n$, where $\lambda_3 \simeq 0.372$, $\lambda_4 \simeq 0.389$, $\lambda_5 \simeq 0.393$ and $\lambda_6 \simeq 0.3935$. These values are slightly sensitive to the tolerance used in the Newton--Raphson algorithm, and to the precision of the calculation. To obtain these orbits numerically, it is useful to note that 
they can be distinguished from the coexisting period-$3$ orbits by the respective values of their maximal Lyapunov exponent. Although this is always negative
for periodic orbits, its value is around $-0.3$ for the attracting period-$3$ orbits, and around $-0.01$ for period-$32$, for example.

By  $\lambda = 0.3937$ a chaotic attractor (i.e.\ one with a positive Lyapunov exponent) is obtained. This chaotic attractor has features which appear visually similar to the type of structure found in the H\'enon attractor \cite{HenonTwoDimensionalMappingStrangeAttractor1976}. Part of the bifurcation diagram is shown in \figref{fig:bifn} for $\alpha=0.08$, where the period-doubling cascade and the subsequent merging crisis (see below) can be seen.  

\begin{figure}
\includegraphics{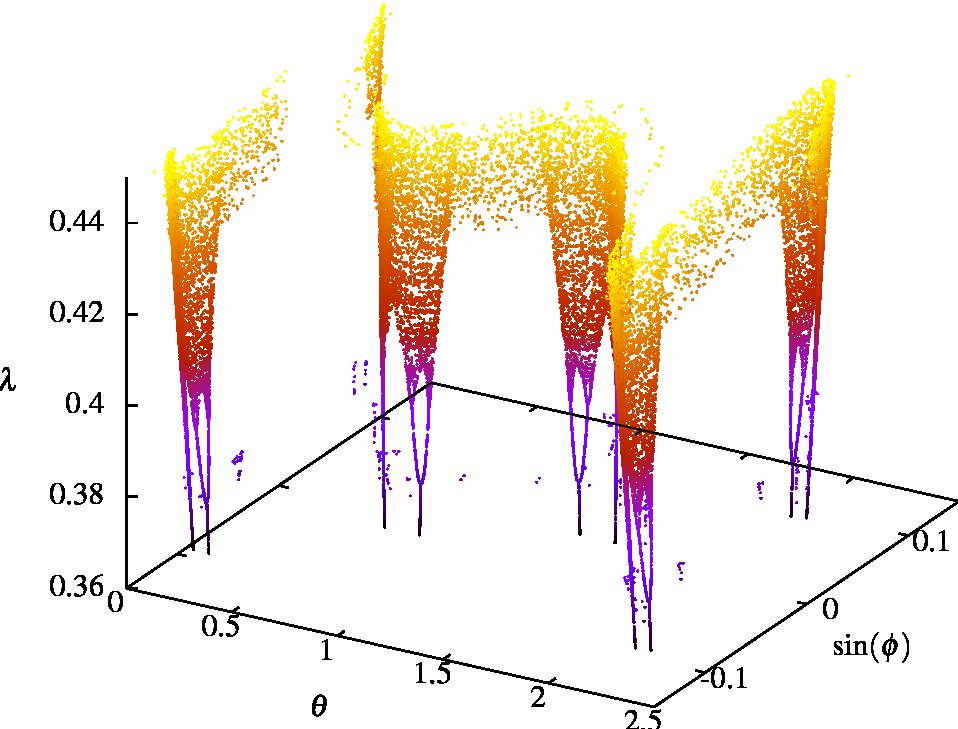}
\caption{Part of the bifurcation diagram of the three-pointed egg with $\alpha = 0.08$,
for a range of $\lambda$ including the periodic-doubling cascade, and the subsequent crisis of the chaotic attractors when three spatially distinct attractors combine into a single one. For each horizontal slice with fixed $\lambda$, attractors in phase space obtained from different initial conditions are shown, indicating the range of possible behaviour for each value of $\lambda$.
The ``noise'' visible for smaller values of $\lambda$ are period-doubling cascades of other orbits, which also reach a crisis when they are absorbed into other attractors.
} 
\label{fig:bifn}
\end{figure}

Also visible on \figref{fig:polar-phase-diag} are what appear to be isolated regions with chaotic attractors. 
The fact that they appear isolated is apparently due to the resolution of the plot.
These chaotic attractors also appear to be created by period-doubling cascades, starting from periodic orbits of higher period, as is visible as what initially appears to be ``noise'' in \figref{fig:bifn}.  In this work we have not attempted to exhaustively characterise the structure of attractors, which are known to exhibit very complicated coexistence phenomena in other models \cite{FeudelGrebogiMap100CoexistingPeriodicAttractorsPRE1996, GallasMultistabilityKickedRotorIJBC2008}.

\subsection{Crises of chaotic attractor}
We continue to fix $\alpha = 0.08$.
Above $\lambda_c \simeq 0.42$, the attractor has behaviour reminiscent of the logistic map, with 
alternating chaotic and periodic windows.
As $\lambda$ increases, the chaotic attractor expands in phase space, and has an appearance somewhat similar to that of the H\'enon attractor. 

Initially there are three separate attractors, each concentrated near one of the diametric period-$2$ orbits. 
However, at around $\lambda = 0.4435$, there is a first crisis \cite[Chapter~8]{OttChaosDynamicalSystemsBook1993}  in the chaotic attractor: the previously separate three symmetric parts of the attractor abruptly join together into a single attractor.

\begin{figure}
\subfigure[$\lambda = 0.66$]{
\includegraphics[scale=0.9]{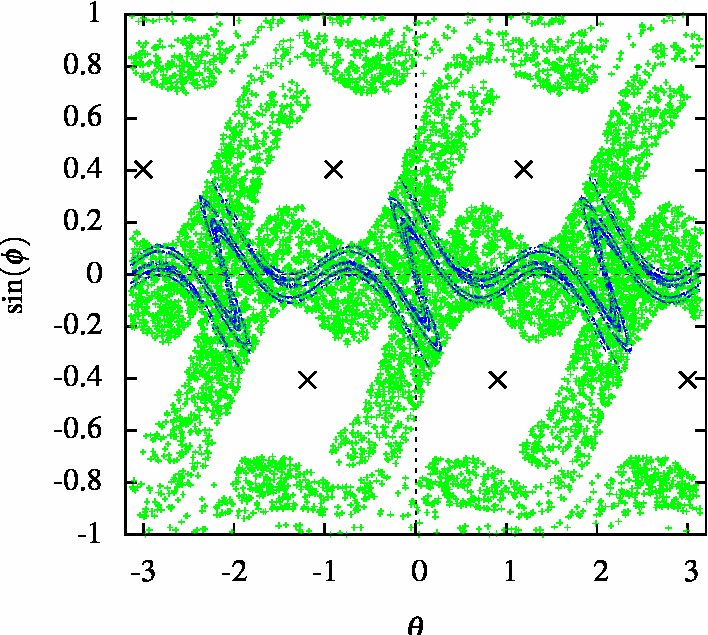}
 \label{fig:crisis}
}
\qquad
\subfigure[$\lambda = 0.7$]{
\includegraphics[scale=0.9]{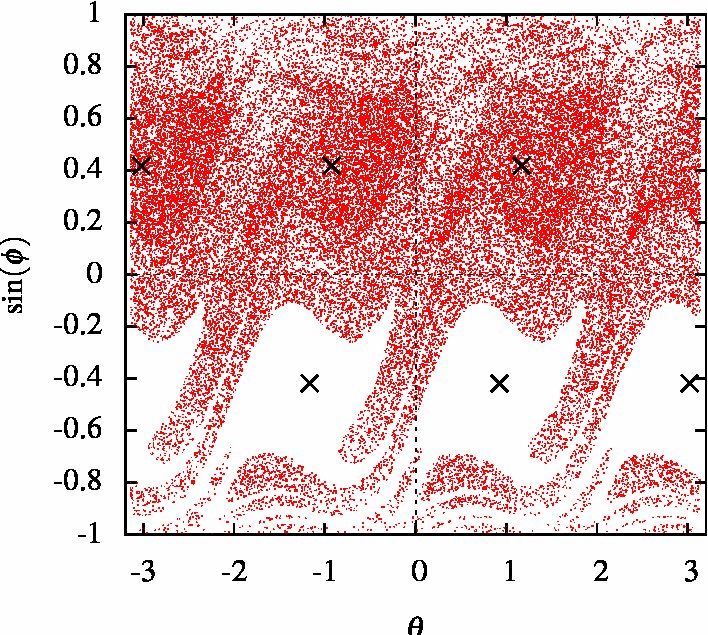}
\label{fig:intermingled}
}
\caption{Basins of attraction in phase space for the three-pointed egg with $\alpha=0.08$ and 
(a) $\lambda = 0.66$; (b) $\lambda = 0.7$. In (a), the chaotic attractor is shown with blue points and its
basin of attraction with green $+$; the period-$3$ attracting orbits are shown as black $\times$, and their basins of attraction in white. The individual basins of the two period-$3$ orbits are disjoint, with one lying above and the other below the chaotic attractor.
In (b), the chaotic attractor has disappeared in a crisis, leaving just two attracting period-$3$ orbits, shown as black $\times$. The basins of each of these orbits are distinguished, one in red dots and the other blank.
The region that in (a) was occupied by the basin of the chaotic attractor is now a region where the basins of the two periodic orbits intermingle, apparently arbitrarily closely in the whole region.
}
\end{figure}

This single attractor then continues to expand in phase space until
another crisis occurs: the chaotic attractor abruptly vanishes
when $\lambda \simeq 0.66$. The reason for this is found in \figref{fig:crisis}, which shows the
basins of attraction of 
the chaotic and period-$3$ attractors 
for $\lambda=0.66$, just below the crisis.
As $\lambda$ increases, both the chaotic attractor and the basins of the periodic
attractors grow. For a slightly larger value of $\lambda$, the chaotic attractor touches the basins of the period-$3$ attractors. After a chaotic transient which resembles the previously-existing attractor, all initial 
conditions are then attracted to the period-$3$ orbits, and the chaotic attractor thus disappears.

In figure~\ref{fig:crisis}, the basins of the two distinct period-$3$ orbits are not distinguished. In fact, they are completely separate: the points close to each orbit are attracted to it, and the basins are divided by that of the chaotic attractor. Once the chaotic attractor has vanished, however, the region that was occupied by its basin is now attracted to one of the two period-$3$ orbits. Their individual basins are shown in 
figure~\ref{fig:intermingled}. We find that their basins  are apparently completely intermingled in this intermediate region.

We remark that the chaotic attractor found in the three-pointed egg is interesting since
it is the only case that we have studied of a billiard table having
a \emph{smooth} boundary (that is, with no points at which the
boundary is not $C^3$) which exhibits a chaotic attractor.

\section{Discussion and conclusions}
 
We have shown that there is a large variety of interesting phenomena in dissipative pinball billiards. Generically, a mixture of periodic and chaotic attractors is found,
with only the former in the ellipse and only the latter in the cardioid.
 
In the cuspless cardioid, we showed that discontinuities in the curvature of the table boundary can lead to the the sudden creation of a chaotic attractor. These discontinuities are also reflected in the existence of a kink structure in the resulting attractor in phase space.

Smooth tables, such as the three-pointed egg, can exhibit complicated sequences of bifurcations, including period-doubling cascades from periodic orbits of different periods. These can lead to chaotic attractors which then collapse by means of crises when they touch a different basin of attraction.

In \cite{MarkarianPujalsSambarinoPinballBilliards}, it was shown that in systems
with dominated splitting there is another possibility in addition to periodic
and chaotic attractors, namely the existence of sets diffeomorphic to a circle
where the dynamics is conjugate to an irrational rotation (for details, see
\cite[Section~2]{MarkarianPujalsSambarinoPinballBilliards}).
We have not observed such behaviour in the models we have studied.
One of the possible reasons for the absence of such behaviour is the fact that
we have only used uniform contraction of the angle at collisions; other types
of collision rule are under study.
This kind of behaviour could, however, perhaps be found on the boundaries which
divide different basins of attraction, which are often found to be smooth
one-dimensional curves, as depicted in several of the figures in this paper.

We found that weakly-dissipative pinball billiards, i.e.\ with $\lambda$ close to $1$, give some information on the existence of elliptic periodic orbits in the elastic case.
Furthermore, we studied two cases, the cardioid and the cuspless cardioid, for which the dynamics is fully chaotic when $\lambda=1$, and in this case weak dissipation leads instead to a chaotic attractor. We suggest that adding weak dissipation could thus
possibly be useful as a means of detecting small elliptic islands in the phase
space of Hamiltonian billiards. An interesting alternative method for doing
this was proposed in
ref.~\cite{TailleurKurchanRarePhysicalLyapunovDynamicsNatPhys2007}.


An open question is whether a similar dynamical portrait in parameter space occurs for similar billiards, for example a $K$-pointed egg given by the polar equation $\rho(\theta) = 1 + \alpha \cos(K \theta)$, for other values of $K>3 \in \mathbb{N}$.

The study of pinball billiard dynamics in other types of billiard table, 
including dispersing billiards, is a work in progress.

\ack
RM and DPS thank the Erwin Schr\"odinger Institute, Vienna, for financial support, and the Unidad Cuernavaca of the Instituto de Matem\'aticas, Universidad Nacional Aut\'onoma de M\'exico, for hospitality. DPS acknowledges support from the PROFIP programme of DGAPA-UNAM.
Support from PAPIIT project IN102307 (DGAPA-UNAM), CONACYT grant no.~58354, PDT project S/C/IF/54/001,  PEDECIBA, Uruguay, and the Academia Mexicana de Ciencias is also acknowledged.

\appendix
\section{Numerical methods}

In this appendix we summarise the numerical methods used to calculate the figures shown in the paper.

\subsection{Trajectories and attractors}

Initial conditions for the simulations are generally chosen uniformly in the interior of the billiard table, with initial velocities chosen with uniform directions $\phi \in [0, 2\pi)$, and are then mapped to the first intersection with the table boundary.  The subsequent dynamics is followed by finding each consecutive intersection with the boundary and implementing the pinball-type collisions. Cartesian coordinates are used whenever possible.

The distance (or equivalently time) $t$ until the next intersection with the table boundary is found using techniques specific to each type of table.  For the ellipse, the intersection is given by solving a quadratic equation which is a generalisation of that used for intersections with a circle \cite{GaspardBook1998}. In the cardioid, a Newton--Raphson step is used to find the first boundary collision, and thereafter a cubic equation is solved, following \cite{BackerDullinSymbDynCardioidBilliardJPA1997}. In the cuspless cardioid, this is adapted so that if the proposed root lies  in the forbidden region then the particle instead collides with the plane.  In the three-pointed egg, a Newton--Raphson method is used to find the intersections.

Once the collision point $\qq$ on the boundary has been found, the pinball billiard dynamics is implemented as follows.
The normal vector $\nn(\qq)$ at the point $\qq$ of collision is calculated based on the equation describing the table boundary, and a tangent vector $\ttt(\qq)$ there is constructed orthogonal to $\nn$.
The incoming angle $\eta$ at the collision is found as $\eta = \cos^{-1} (\vv \cdot \nn)$, and the outgoing angle is then given by $\phi = \lambda \eta$. The outgoing velocity in Cartesian coordinates is finally given by $\vv' = \nn \cos \phi  +  \ttt d \sin \phi$, where $d = \pm 1$ ensures that the component of the outgoing velocity is in the same direction with respect to the tangent vector as the incoming velocity was, as determined by the sign of $\vv \cdot \ttt$.
Note that in the case of elastic collisions, the dynamics may be implemented completely in Cartesian coordinates \cite{GaspardBook1998}; however, the definition of pinball billiard dynamics makes the use of trigonometric functions unavoidable.

To find attractors of the dynamics, the first several thousand collisions in the trajectory emanating from each initial collision are discarded.  Assuming that transients have already decayed by that point, and that the dynamics is transitive on the attractor, the result is an approximate picture of the attractor.

\subsection{Lyapunov exponents}

The Lyapunov exponents of an orbit of the pinball billiard map $T_\lambda$ are calculated using the algorithm of Dellago et al.\ \cite{DellagoPoschLyapInstabHardDisksEqmNoneqPRE1996}. This involves the derivative \eqref{eq:derivative_of_map} of the map of the pinball billiard transformation.
For completeness we summarise the main steps of the procedure.

Given an initial point $\xx$ in the two-dimensional phase space, several thousand collisions in its trajectory are discarded to remove transients. Then the trajectory is extended to the next collision, and two orthonormal two-dimensional vectors $\vv_1$ and $\vv_2$, which are thought of as tangent vectors at the point $\xx$, are manipulated in parallel. Each tangent vector is acted on by the derivative matrix $D_{\xx} T_\lambda$, to give two new vectors $\vv_1'$ and $\vv_2'$, which are then reorthogonalised using a Gram--Schmidt procedure.  This process is repeated at each collision.  The Lyapunov exponents $\nu_1$ and $\nu_2$ are given by the average growth rates of the normalising factors.  In particular, their sum $\nu_1 + \nu_2$ gives the logarithm of the mean rate of volume contraction of the dynamics on the attractor in phase space.

 
\subsection{Basins of attraction}
To determine basins of attraction in phase space, 
initial conditions are chosen uniformly in phase space as above (with the initial boundary collision being an elastic one, to allow all initial directions to be produced).
For each initial condition $(s, \phi)$, the two Lyapunov exponents are calculated. If the maximal Lyapunov exponent is positive, then the trajectory converges to a chaotic attractor; if it is negative or $0$ then it converges to a periodic attractor.  A symbol is plotted at $(s, \phi)$ according to the type of attractor found there.
For the correspondence between negative Lyapunov exponents and periodic orbits, and between positive Lyapunov exponents and chaotic behaviour on strange attractors, we refer to
\cite[Section III]{EckmannRuelleErgTheoryChaosStrangeAttractorsRMP1985}.

\section*{References}

\def\cprime{$'$}

\end{document}